\documentclass[a4paper,11pt]{amsart}
\usepackage{amsfonts,amsmath,amsthm,amssymb}
\usepackage{enumitem}
\usepackage{mathrsfs}
\usepackage[margin=1in]{geometry}
\usepackage[colorlinks=true,linkcolor=blue,filecolor=blue,citecolor=red]{hyperref}

\usepackage{cleveref}
\usepackage{amsrefs}

\theoremstyle{plain}
\newtheorem{theorem}{Theorem}[section]
\newtheorem{corollary}[theorem]{Corollary}
\newtheorem{lemma}[theorem]{Lemma}

\newtheorem{proposition}[theorem]{Proposition}

\newtheorem{conjecture}[theorem]{Conjecture}
\theoremstyle{definition}

\newtheorem{example}[theorem]{Example}

\begin{document}
\title[Ramanujan's continued fractions of order $10$ as modular functions]{Ramanujan's continued fractions of order $10$ as modular functions}
\author[Victor Manuel Aricheta, Russelle Guadalupe]{Victor Manuel Aricheta, Russelle Guadalupe}
\address{Institute of Mathematics, University of the Philippines, Diliman\\
Quezon City 1101, Philippines}
\email{vmaricheta@math.upd.edu.ph, rguadalupe@math.upd.edu.ph}

\renewcommand{\thefootnote}{}



\renewcommand{\thefootnote}{\arabic{footnote}}
\setcounter{footnote}{0}

\begin{abstract}
We explore the modularity of the continued fractions $I(\tau), J(\tau), T_1(\tau), T_2(\tau)$ and $U(\tau)=I(\tau)/J(\tau)$ of order $10$, where $I(\tau)$ and $J(\tau)$ are introduced by Rajkhowa and Saikia, which are special cases of certain identities of Ramanujan. In particular, we show that these fractions can be expressed in terms of an $\eta$-quotient $g(\tau)$ that generates the field of all modular functions on the congruence subgroup $\Gamma_0(10)$. Consequently, we prove that modular equations for $g(\tau)$ and $U(\tau)$ exist at any level and derive these equations of prime levels $p\leq 11$. We also show that the continued fractions of order $10$ can be explicitly evaluated using a singular value of $g(\tau)$, which under certain conditions, generates the Hilbert class field of an imaginary quadratic field. We employ the methods of Lee and Park to establish our results.	
\end{abstract}

\maketitle

\section{Introduction}\label{sec1}

\subsection{Background on $q$-continued fractions and their modular properties}
Throughout this paper, let $q:= e^{2\pi i\tau}$ for an element $\tau$ in the complex upper half-plane $\mathbb{H}$ and $(a;q)_{\infty} :=\prod_{n=0}^{\infty}(1-aq^n)$ for $a\in\mathbb{C}$. 
We define the Rogers-Ramanujan continued fraction by
\begin{align*}
r(\tau) = \dfrac{q^{1/5}}{1+\dfrac{q}{1+\dfrac{q^2}{1+\dfrac{q^3}{1+\cdots}}}} = q^{1/5}\dfrac{(q;q^5)_{\infty}(q^4;q^5)_{\infty}}{(q^2;q^5)_{\infty}(q^3;q^5)_{\infty}}.
\end{align*}
In his first two letters \cite[p. 21 -- 30, p. 53 -- 62]{ramlet} to Hardy in 1913, Ramanujan claimed the following exact values
\begin{align*}
r(i) &= \sqrt{\dfrac{5+\sqrt{5}}{2}}-\dfrac{\sqrt{5}+1}{2},\\
r(\sqrt{-5}) &= \dfrac{\sqrt{5}}{1+(5^{3/4}(\frac{\sqrt{5}-1}{2})^{5/2}-1)^{1/5}}-\dfrac{\sqrt{5}+1}{2},
\end{align*}
and asserted that the values of $r(\sqrt{-n}/2)$ \textquotedblleft can be exactly found\textquotedblright\,for any positive rational number $n$. Ramanujan's observations were later proven correct by Berndt, Chan and Zhang \cite{zhchbr}, who found a systematic way of evaluating the exact values of $r(\sqrt{-n}/2)$ for positive rational numbers $n$ by using Ramanujan's eta-function identities. Gee and Honsbeek \cite{geehons} extended their results by proving that $r(\tau)$, viewed as a modular function, generates the field of all modular functions on the congruence subgroup $\Gamma(5)$. Consequently, they showed that $r(\tau)$ generates the ray class field modulo $5$ over an imaginary quadratic field $K$ for some $\tau\in K \cap \mathbb{H}$ and explicitly computed new values of $r(\tau)$ as nested radicals. After the modularity and the arithmetic of $r(\tau)$ were established, other $q$-continued fractions were studied by several authors \cite{chokoop,chokoop2,duke,leeparka,leeparkc,leeparke} for their analogous properties, with the primary goal of evaluating these fractions at certain imaginary quadratic points and constructing abelian extensions over suitable imaginary quadratic fields.

Aside from the fraction $r(\tau)$, Ramanujan recorded several identities involving $q$-continued fractions in his notebooks; one of which is given by \cite[p. 24, Entry 12]{bcberndt}
\begin{align}\label{eq1}
\dfrac{(a^2q^3;q^4)_\infty(b^2q^3;q^4)_\infty}{(a^2q;q^4)_\infty(b^2q;q^4)_\infty}=\dfrac{1}{1-ab+\dfrac{(a-bq)(b-aq)}{(1-ab)(q^2+1)+\dfrac{(a-bq^3)(b-aq^3)}{(1-ab)(q^4+1)+\cdots}}},
\end{align} 
valid for all complex numbers $a$ and $b$ with $|ab| < 1$, or those satisfying $a=b^{2m+1}$ for some integer $m$. Park \cite{park1,park2} studied the modularity of the continued fractions of levels $16$ and $18$, which are particular cases of (\ref{eq1}), and showed that these can be expressed in terms of the generators for the field of all modular functions on $\Gamma_0(16)$ and $\Gamma_0(18)$, respectively. 

On the other hand, Ramanujan also recorded the following identity \cite[p. 21, Entry 11]{bcberndt}
\begin{align}
R(a,b;q)&:=\dfrac{(-a;q)_{\infty}(b;q)_{\infty}-(a;q)_{\infty}(-b;q)_{\infty}}{(-a;q)_{\infty}(b;q)_{\infty}+(a;q)_{\infty}(-b;q)_{\infty}}\nonumber\\
&=\dfrac{a-b}{1-q+\dfrac{(a-bq)(aq-b)}{1-q^3+\dfrac{q(a-bq^2)(aq^2-b)}{1-q^5+\cdots}}},\label{eq2}
\end{align}
where $a, b\in\mathbb{C}$ with $|ab|<1$. Lee and Park \cite{leeparkd} showed that the continued fractions of order $4$ given by
\begin{align*}
T(\tau) &:=R(q,0;q^2)=\dfrac{q}{1-q^2+\dfrac{q^4}{1-q^6+\dfrac{q^8}{1-q^{10}+\cdots}}}\\
&=\dfrac{(-q;q^2)_{\infty}-(q;q^2)_{\infty}}{(-q;q^2)_{\infty}+(q;q^2)_{\infty}}
\end{align*}
and
\begin{align*}
W(\tau)&:=\dfrac{1}{2}R(q,-q;q^2)=\dfrac{q}{1-q^2+\dfrac{q^2(1+q^2)^2}{1-q^6+\dfrac{q^4(1+q^4)^2}{1-q^{10}+\cdots}}}\\
&=\dfrac{1}{2}\dfrac{(-q;q^2)_{\infty}^2-(q;q^2)_{\infty}^2}{(-q;q^2)_{\infty}^2+(q;q^2)_{\infty}^2},
\end{align*}
which were first defined by Saikia \cite{saikia1,saikia2} and are particular cases of (\ref{eq2}), are modular functions on $\Gamma_0(32)$ and $\Gamma_0(16)$, respectively. They also proved that $W(\tau/4)$ generates the ray class field modulo $4$ over an imaginary quadratic field $K$ for some $\tau\in K \cap \mathbb{H}$ and that $T(\tau)$ is an algebraic unit.

Recently, the second author \cite{guad} applied the methods of Lee and Park \cite{leeparkb,leeparkd} to show that the continued fraction of order $6$ given by
\begin{align*}
H(\tau):=R(q,-q^2;q^3)&=\dfrac{q(1+q)}{1-q^3+\dfrac{q^3(1+q^2)(1+q^4)}{1-q^9+\dfrac{q^6(1+q^5)(1+q^7)}{1-q^{15}+\cdots}}}\\
&=\dfrac{(-q;q^3)_{\infty}(-q^2;q^3)_{\infty}-(q;q^3)_{\infty}(q^2;q^3)_{\infty}}{(-q;q^3)_{\infty}(-q^2;q^3)_{\infty}+(q;q^3)_{\infty}(q^2;q^3)_{\infty}},
\end{align*}
which is another special case of (\ref{eq2}), is a modular function on $\Gamma_0(24)$ that can be written in terms of a generator $s(\tau)$ for the field of all modular functions on $\Gamma_0(12)$.
He also proved that $s(\tau/2)$ generates the ray class field modulo $6$ over an imaginary quadratic field $K$ for some $\tau\in K \cap \mathbb{H}$ and that $H(\tau)$ is an algebraic unit.

\subsection{Objectives of this paper}

In this paper, we study the modular properties of the continued fractions of order $10$ given by $U(\tau):=I(\tau)/J(\tau)$, where 
\begin{align*}
I(\tau) &:= q^{3/4}\dfrac{(q;q^{10})_\infty(q^9;q^{10})_\infty}{(q^4;q^{10})_\infty(q^6;q^{10})_\infty}\\
&= \dfrac{q^{3/4}(1-q)}{1-q^{5/2}+\dfrac{q^{5/2}(1-q^{3/2})(1-q^{7/2})}{(1-q^{5/2})(1+q^5)+\dfrac{q^{5/2}(1-q^{13/2})(1-q^{17/2})}{(1-q^{5/2})(1+q^{10})+\cdots}}}
\end{align*}
and 
\begin{align*}
J(\tau) &:= q^{1/4}\dfrac{(q^2;q^{10})_\infty(q^8;q^{10})_\infty}{(q^3;q^{10})_\infty(q^7;q^{10})_\infty}\\
&= \dfrac{q^{1/4}(1-q^2)}{1-q^{5/2}+\dfrac{q^{5/2}(1-q^{1/2})(1-q^{9/2})}{(1-q^{5/2})(1+q^5)+\dfrac{q^{5/2}(1-q^{11/2})(1-q^{19/2})}{(1-q^{5/2})(1+q^{10})+\cdots}}}
\end{align*}
are two new special cases of (\ref{eq1}) recently introduced by Rajkhowa and Saikia \cite{rajsai}, and
\begin{align*}
T_1(\tau) &:=R(q,-q^4;q^5)=\dfrac{q(1+q^3)}{1-q^5+\dfrac{q^5(1+q^2)(1+q^8)}{1-q^{15}+\dfrac{q^{10}(1+q^7)(1+q^{13})}{1-q^{25}+\cdots}}}\\
&=\dfrac{(-q;q^5)_{\infty}(-q^4;q^5)_{\infty}-(q;q^5)_{\infty}(q^4;q^5)_{\infty}}{(-q;q^5)_{\infty}(-q^4;q^5)_{\infty}+(q;q^5)_{\infty}(q^4;q^5)_{\infty}},\\
T_2(\tau)&:=R(q^2,-q^3;q^5)=\dfrac{q^2(1+q)}{1-q^5+\dfrac{q^5(1+q^4)(1+q^6)}{1-q^{15}+\dfrac{q^{10}(1+q^9)(1+q^{11})}{1-q^{25}+\cdots}}}\\
&=\dfrac{(-q^2;q^5)_{\infty}(-q^3;q^5)_{\infty}-(q^2;q^5)_{\infty}(q^3;q^5)_{\infty}}{(-q^2;q^5)_{\infty}(-q^3;q^5)_{\infty}+(q^2;q^5)_{\infty}(q^3;q^5)_{\infty}},
\end{align*}
which are two new special cases of (\ref{eq2}). More precisely, we show that these fractions can be expressed in terms of an $\eta$-quotient defined by 
\begin{equation*}
g(\tau):=\dfrac{\eta^2(2\tau)\eta^4(5\tau)}{\eta^4(\tau)\eta^2(10\tau)},
\end{equation*}
where $\eta(\tau):=q^{1/24}(q;q)_\infty$ is the Dedekind eta function, that generates the field of all modular functions on $\Gamma_0(10)$. We then extend the results of \cite{rajsai} by using the modular properties of $g(\tau)$ and $U(\tau)$ to derive modular equations for these functions at any level, which can be used to find the exact values of $I(\tau), J(\tau), T_1(\tau)$ and $T_2(\tau)$ for imaginary quadratic points $\tau$. We employ the methods of Lee and Park \cite{leeparkb,leeparkd} and certain properties of $\eta$-quotients and generalized $\eta$-quotients defined by Yang \cite{yang} to present our results.

The paper is organized as follows. Section \ref{sec2} gives a brief background about modular functions on $\Gamma_0(N)$ and $\Gamma_1(N)$, with an emphasis on $\eta$-quotients and generalized $\eta$-quotients. Section \ref{sec3} deals with the modular properties of the functions $I(\tau), J(\tau), T_1(\tau), T_2(\tau), U(\tau)$ and $g(\tau)$. We show in Section \ref{sec4} that for positive integers $n$, the modular equation for $g(\tau)$ of level $n$ exist, which serves as an affine model for the modular curve $X_0(10n)$ over $\mathbb{Q}$, and provide modular equations for $g(\tau)$ and $U(\tau)$ of prime levels $p\leq 11$. We prove in Section \ref{sec5} that if $K$ is an imaginary quadratic field and $\tau\in K\cap\mathbb{H}$ is a root of a primitive polynomial $aX^2+bX+c\in\mathbb{Z}[X]$ satisfying certain conditions, then $g(\tau)$ generates the Hilbert class field of $K$. We also show that $1/g(\tau), 1/U(\tau), 1/I(\tau)$ and $J(\tau)$ are algebraic integers, and $T_1(\tau)$ and $T_2(\tau)$ are algebraic units for $\tau\in K\cap \mathbb{H}$. We finally provide an example of finding the minimal polynomial of $1/g(\tau)$ over $K$, using Shimura reciprocity law due to Cho and Koo \cite{chokoo} and Gee \cite{gee}. We then use the exact value of $g(\tau)$ to find $U(\tau), I(\tau), J(\tau), T_1(\tau)$ and $T_2(\tau)$. We have performed most of our computations via \textit{Mathematica}.

\section{Modular functions on $\Gamma_0(N)$ and $\Gamma_1(N)$}\label{sec2}

In this section, we give a basic background about modular functions on congruence subgroups, which are subgroups of $\mbox{SL}_2(\mathbb{Z})$ containing 
\begin{equation*}
\Gamma(N) := \left\lbrace\begin{bmatrix}
a & b\\ c& d
\end{bmatrix} \in \mbox{SL}_2(\mathbb{Z}) : a-1\equiv b\equiv c\equiv d-1\equiv 0\pmod N\right\rbrace
\end{equation*}
for some integer $N\geq 1$. Examples of congruence subgroups we focus on in this paper are
\[\begin{aligned}
\Gamma_1(N) &:= \left\lbrace\begin{bmatrix}
a & b\\ c& d
\end{bmatrix} \in \mbox{SL}_2(\mathbb{Z}) : a-1\equiv c\equiv d-1\equiv 0\pmod N\right\rbrace,\\
\Gamma_0(N) &:= \left\lbrace\begin{bmatrix}
a & b\\ c& d
\end{bmatrix} \in \mbox{SL}_2(\mathbb{Z}) : c\equiv 0\pmod N\right\rbrace.
\end{aligned}\]
An element of a given congruence subgroup $\Gamma$ acts on the extended upper half-plane $\mathbb{H}^\ast := \mathbb{H}\cup \mathbb{Q}\cup \{\infty\}$ by a linear fractional transformation. We define the cusps of $\Gamma$ as the equivalence classes of $\mathbb{Q}\cup \{\infty\}$ under this action, and it is known that there are finitely many inequivalent cusps of $\Gamma$. We define a modular function on $\Gamma$ as a function $f:\mathbb{H}\rightarrow\mathbb{C}$ such that $f(\tau)$ is meromorphic on $\mathbb{H}$, $f(\gamma\tau)=f(\tau)$ for all $\gamma\in \Gamma$, and $f(\tau)$ is meromorphic at all cusps of $\Gamma$. The last condition means that for every cusp $r$ of $\Gamma$ and an element $\gamma\in\mbox{SL}_2(\mathbb{Z})$ with $\gamma(\infty)=r$, the $q$-expansion of $f\circ \gamma$ is given by $(f\circ\gamma)(\tau) = \sum_{n\geq n_0} a_nq^{n/h}$ for some integers $h$ and $n_0$ with $a_{n_0}\neq 0$. We call $n_0$ the order of $f(\tau)$ at $r$, denoted by $\mbox{ord}_r f(\tau)$, and we say that $f(\tau)$ has a zero (resp., a pole) at $r$ if $\mbox{ord}_r f(\tau)$ is positive (resp., negative). We also call $h$ the width of $r$, which is the smallest positive integer such that $\gamma[\begin{smallmatrix}
1 & h\\0 & 1
\end{smallmatrix}]\gamma^{-1}\in \pm\Gamma$. We note that $h$ is independent on the choice of $\gamma$.

Let $A_0(\Gamma)$ be the field of all modular functions on $\Gamma$ and $A_0(\Gamma)_{\mathbb{Q}}$ be its subfield consisting of all modular functions on $\Gamma$ with rational $q$-expansion. We identify $A_0(\Gamma)$ with the field of all meromorphic functions on the modular curve $X(\Gamma):=\Gamma\backslash \mathbb{H}^\ast$. If $f(\tau)\in A_0(\Gamma)$ has no zeros nor poles on $\mathbb{H}$, then the degree $[A_0(\Gamma):\mathbb{C}(f(\tau))]$ is the total degree of poles of $f(\tau)$ given by $-\sum_r \mbox{ord}_r f(\tau)$, where the sum ranges over the inequivalent cusps $r$ of $\Gamma$ for which $f(\tau)$ has a pole at $r$ (see \cite[Proposition 2.11]{shimura}). The following lemma gives a complete set of all inequivalent cusps of $\Gamma_0(N)$ with their respective widths.

\begin{lemma}[\cite{chokoop}]\label{lem21}
Let $a,a',c$ and $c'$ be integers with $\gcd(a,c)=\gcd(a',c')=1$. Let $S_{\Gamma_0(N)}$ be the set of all inequivalent cusps of $\Gamma_0(N)$. We denote $\pm1/0$ as $\infty$. Then 
\begin{enumerate}
\item $a/c$ and $a'/c'$ are equivalent over $\Gamma_0(N)$ if and only if there exist an integer $n$ and an element $\overline{s}\in (\mathbb{Z}/N\mathbb{Z})^\times$ such that $(a',c')\equiv(\overline{s}^{-1}a+nc,\overline{s}c)\pmod N$,
\item we have 
\begin{align*}
	S_{\Gamma_0(N)} = \{a_{c,j}/c\in \mathbb{Q} : 0 < c \mid N, 0 < a_{c,j}\leq N, \gcd(a_{c,j},N)=1,\\ a_{c,j}= a_{c,j'}\stackrel{\text{def}}{\iff} a_{c,j}\equiv a_{c,j'}\pmod{\gcd(c,N/c)}\}, and
\end{align*}
\item the width of a cusp $a/c\in S_{\Gamma_0(N)}$ is $N/\gcd(c^2,N)$.
\end{enumerate}
\end{lemma}

We define an $\eta$-quotient by a function of the form
\begin{equation*}
f(\tau) = \prod_{\delta\mid N} \eta(\delta\tau)^{r_{\delta}}
\end{equation*}
for some indexed set $\{r_\delta\in\mathbb{Z} : \delta\mid N\}$. The next two lemmas, due to Newman and Ligozat, provide sufficient conditions for the modularity of an $\eta$-quotient.

\begin{lemma}[\cite{newman1, newman2,ono}]\label{lem22}
Let $f(\tau) = \prod_{\delta\mid N} \eta(\delta\tau)^{r_{\delta}}$ be an eta-quotient with $k= \frac{1}{2}\sum_{\delta\mid N} r_{\delta}\in\mathbb{Z}$ such that
\begin{equation*}
\sum_{\delta\mid N} \delta r_{\delta} \equiv 0\pmod{24}\quad\text{ and }\quad\sum_{\delta\mid N} \dfrac{N}{\delta}r_{\delta} \equiv 0\pmod {24}.
\end{equation*}
Then for all $[\begin{smallmatrix}
a & b\\ c & d
\end{smallmatrix}]\in \Gamma_0(N)$,
\begin{equation*}
f\left(\dfrac{a\tau+b}{c\tau+d}\right) = \left(\dfrac{(-1)^k\prod_{\delta\mid N} \delta^{r_\delta}}{d}\right)(c\tau+d)^kf(\tau),
\end{equation*}
where $(\frac{\cdot}{d})$ is the Legendre symbol.
\end{lemma} 

\begin{lemma}[\cite{ligo,ono}]\label{lem23}
Let $c, d$ and $N$ be positive integers with $d\mid N$ and $\gcd(c,d)=1$ and let $f(\tau) = \prod_{\delta\mid N} \eta(\delta\tau)^{r_{\delta}}$ be an $\eta$-quotient satisfying the conditions of Lemma \ref{lem22}. Then the order of vanishing of $f(\tau)$ at the cusp $c/d$ is 
\begin{align*}
\dfrac{N}{24d\gcd(d,\frac{N}{d})}\sum_{\delta\mid N} \gcd(d,\delta)^2\cdot\dfrac{r_{\delta}}{\delta}.
\end{align*}
\end{lemma}

We next define, following Yang \cite{yang}, generalized $\eta$-quotients with analogous properties on $\Gamma_1(N)$. These are functions of the form
\begin{equation*}
f(\tau) = \prod_{1\leq g\leq \lfloor N/2\rfloor}\eta_{N,g}(\tau)^{r_g}
\end{equation*}
for some indexed set $\{r_g \in\mathbb{Z}: 1\leq g\leq \lfloor N/2\rfloor\}$ of integers, where 
\begin{equation*}
\eta_{N, g}(\tau) = q^{NB(g/N)/2}\prod_{m=1}^\infty(1-q^{N(m-1)+g})(1-q^{Nm-g})
\end{equation*}
is the generalized Dedekind eta function with $B_2(t) := t^2-t+1/6$. The following results give the transformation formula of $\eta_{N, g}$ on $\Gamma_0(N)$ and conditions for the modularity of a generalized $\eta$-quotient on $\Gamma_1(N)$.

\begin{lemma}[\cite{yang}]\label{lem24}
The function $\eta_{N,g}(\tau)$ satisfies $\eta_{N,g+N}(\tau)=\eta_{N,-g}(\tau)=-\eta_{N,g}(\tau)$. Moreover, let $\gamma = [\begin{smallmatrix}
a & b\\ cN & d
\end{smallmatrix}]\in \Gamma_0(N)$ with $c\neq 0$. Then we have 
\begin{align*}
\eta_{N,g}(\gamma\tau) = \varepsilon(a,bN,c,d)e^{\pi i(g^2ab/N - gb)}\eta_{N,ag}(\tau),
\end{align*}
where 
\begin{align*}
\varepsilon(a,b,c,d) = \begin{cases}
	e^{\pi i(bd(1-c^2)+c(a+d-3))/6}, & \text{ if }c\equiv 1\pmod{2},\\
	-ie^{\pi i(ac(1-d^2)+d(b-c+3))/6}, & \text{ if }c\equiv 0\pmod{2}.
\end{cases}
\end{align*}
\end{lemma}

\begin{lemma}[\cite{yang}]\label{lem25}
Suppose $f(\tau)=\prod_{1\leq g\leq \lfloor N/2\rfloor}\eta_{N,g}(\tau)^{r_g}$ is a generalized $\eta$-quotient such that
\begin{align*}
\sum_{1\leq g\leq \lfloor N/2\rfloor}r_g\equiv 0\pmod{12}\quad\text{ and }\quad\sum_{1\leq g\leq \lfloor N/2\rfloor}gr_g\equiv 0\pmod{2}.
\end{align*}
Then $f(\tau)$ is a modular function on $\Gamma(N)$. Moreover, if in addition, $f(\tau)$ satisfies
\begin{align*}
\sum_{1\leq g\leq \lfloor N/2\rfloor}g^2r_g\equiv 0\pmod{2N}
\end{align*}
then $f(\tau)$ is a modular function on $\Gamma_1(N)$. 
\end{lemma}

\begin{lemma}[\cite{yang}]\label{lem26}
Let $N$ be a positive integer and $\gamma=[\begin{smallmatrix}
a & b\\ c& d
\end{smallmatrix}]\in\mbox{SL}_2(\mathbb{Z})$. Then the first term of the $q$-expansion of $\eta_{N, g}(\gamma\tau)$ is $\varepsilon q^{\delta}$, where $|\varepsilon|=1$ and 
\begin{equation*}
\delta = \dfrac{\gcd(c,N)^2}{2N}P_2\left(\dfrac{ag}{\gcd(c,N)}\right)
\end{equation*}
with $P_2(t):=B_2(\{t\})$ and $\{t\}$ is the fractional part of $t$.
\end{lemma}

Using Lemmata \ref{lem21} and \ref{lem26}, we can compute the order of a generalized $\eta$-quotient at a cusp of $\Gamma_0(N)$.

\section{Modularity of continued fractions of order $10$}\label{sec3}

We establish in this section the modularity of $g(\tau), I(\tau), J(\tau), T_1(\tau), T_2(\tau)$ and $U(\tau)$, which can be used to derive modular equations for $g(\tau)$ and $U(\tau)$ at any level, and evaluate explicitly the singular values of $I(\tau), J(\tau), T_1(\tau)$ and $T_2(\tau)$.  We first recall the Ramanujan's general theta function 
\begin{equation*}
f(a,b) = \sum_{n=-\infty}^{\infty} a^{n(n+1)/2}b^{n(n-1)/2}
\end{equation*}
valid for complex numbers $a$ and $b$ with $|ab| < 1$. We note that $f(a,b)$ satisfies the Jacobi triple product identity \cite[p. 35, Entry 19]{bcberndt}
\begin{equation*}
f(a,b) = (-a;ab)_\infty(-b;ab)_\infty(ab;ab)_\infty.
\end{equation*}
We next define $f(-q):=f(-q,-q^2)=(q;q)_\infty$ and $\chi(q):=(-q;q^2)_\infty$ and present the following identities, as shown below. 

\begin{lemma}\label{lem31}
We have 
\begin{align*}
f(-q,-q^4)f(-q^2,-q^3) &= f(-q)f(-q^5),\\
f(q,q^9)f(q^3,q^7) &= \chi(q)f(-q^5)f(-q^{20}).
\end{align*}
\end{lemma}

\begin{proof}
See \cite[p. 258, Entry 9(vii)]{bcberndt}. We remark that the first identity also follows from the Jacobi triple product identity.
\end{proof}

We use Lemma \ref{lem31} to show that $U(\tau)$ can be written as an $\eta$-quotient.

\begin{lemma}\label{lem32}
We have 
\begin{align*}
U(\tau) = \dfrac{\eta(\tau)\eta(10\tau)^2}{\eta(2\tau)^2\eta(5\tau)}.
\end{align*}
\end{lemma}

\begin{proof}
Using the definition of $U(\tau)$ and Lemma \ref{lem31}, we have
\begin{align*}
U(\tau) = q^{1/2}\dfrac{f(-q,-q^9)f(-q^3,-q^7)}{f(-q^2,-q^8)f(-q^4,-q^6)}= q^{1/2}\dfrac{\chi(-q)f(q^5)f(-q^{20})}{f(-q^2)f(-q^{10})}.
\end{align*}
Since 
\begin{align*}
\chi(-q) = \dfrac{f(-q)}{f(-q^2)}\qquad\text{ and }\qquad f(q) = \dfrac{f^3(-q^2)}{f(-q)f(-q^4)},
\end{align*}
we see that
\begin{align*}
U(\tau) = q^{1/2}\dfrac{f(-q)}{f(-q^2)}\cdot\dfrac{f^3(-q^{10})}{f(-q^5)f(-q^{20})}\cdot\dfrac{f(-q^{20})}{f(-q^{10})}=\dfrac{\eta(\tau)\eta(10\tau)^2}{\eta(2\tau)^2\eta(5\tau)},
\end{align*}
where the last equality follows from $f(-q) = q^{-1/24}\eta(\tau)$.
\end{proof}

We next observe that $I(\tau)$ and $J(\tau)$ are generalized $\eta$-quotients given by $I(\tau) = \eta_{10,1}(\tau)\eta_{10,4}^{-1}(\tau)$ and $J(\tau) = \eta_{10,2}(\tau)\eta_{10,3}^{-1}(\tau)$.

\begin{theorem}\label{thm33}
The functions $I^4(\tau)$ and $J^4(\tau)$ are modular on $\Gamma_1(10)$. On the other hand, the functions $I(4\tau)$ and $J(4\tau)$ are modular on $\Gamma_1(40)$.
\end{theorem} 

\begin{proof}
Observe that the functions $I^4(\tau) = \eta_{10,1}^4(\tau)\eta_{10,4}^{-4}(\tau)$ and $J^4(\tau) = \eta_{10,2}^4(\tau)\eta_{10,3}^{-4}(\tau)$ satisfy the conditions of Lemma \ref{lem25} with $N=10$, so these are modular functions on $\Gamma_1(10)$ by the same lemma. On the other hand, writing $I(4\tau) = \eta_{40,4}(\tau)\eta_{40,16}^{-1}(\tau)$ and $J(4\tau) = \eta_{40,8}(\tau)\eta_{40,12}^{-1}(\tau)$, we see that $I(4\tau)$ and $J(4\tau)$ satisfy the conditions of Lemma \ref{lem25} with $N=40$, so these are modular functions on $\Gamma_1(40)$ by the same lemma again.
\end{proof}

We now define 
\begin{align}
g_1(\tau) &:= \dfrac{1+T_1(\tau)}{1-T_1(\tau)} =  \dfrac{(-q;q^5)_{\infty}(-q^4;q^5)_{\infty}}{(q;q^5)_{\infty}(q^4;q^5)_{\infty}}=\dfrac{\eta_{10,2}(\tau)}{\eta_{5,1}^2(\tau)},\label{eq3}\\
g_2(\tau) &:= \dfrac{1+T_2(\tau)}{1-T_2(\tau)} =  \dfrac{(-q^2;q^5)_{\infty}(-q^3;q^5)_{\infty}}{(q^2;q^5)_{\infty}(q^3;q^5)_{\infty}}=\dfrac{\eta_{10,4}(\tau)}{\eta_{5,2}^2(\tau)},\label{eq4}
\end{align}
and note that
\begin{align*}
g_1^2(\tau)g_2^2(\tau) &= \dfrac{(-q;q^5)^2_{\infty}(-q^2;q^5)^2_{\infty}(-q^3;q^5)^2_{\infty}(-q^4;q^5)^2_{\infty}}{(q;q^5)^2_{\infty}(q^2;q^5)^2_{\infty}(q^3;q^5)^2_{\infty}(q^4;q^5)^2_{\infty}}\\
&=\dfrac{(-q;q)^2}{(-q^5;q^5)^2}\cdot\dfrac{(q^5;q^5)^2}{(q;q)^2} = \dfrac{(q^2;q^2)^2}{(q^{10};q^{10})^2}\cdot\dfrac{(q^5;q^5)^4}{(q;q)^4}\\
&=\dfrac{\eta^2(2\tau)\eta^4(5\tau)}{\eta^4(\tau)\eta^2(10\tau)} = g(\tau).
\end{align*}

\begin{theorem}\label{thm34}
The functions $g_1^4(\tau)$ and $g_2^4(\tau)$ are modular on $\Gamma_1(10)$. On the other hand, the functions $T_1(\tau)$ and $T_2(\tau)$ are modular on $\Gamma_1(40)$.
\end{theorem} 

\begin{proof}
We find that the functions $g_1^4(\tau)=\eta_{10,1}^{-8}(\tau)\eta_{10,2}^4(\tau)\eta_{10,4}^{-8}(\tau)$ and $g_2^4(\tau)=\eta_{10,2}^{-8}(\tau)\eta_{10,4}^4(\tau)\eta_{10,3}^{-8}(\tau)$ satisfy the conditions of Lemma \ref{lem25} with $N=10$, so these are modular functions on $\Gamma_1(10)$ by the same lemma. On the other hand, we observe that
\begin{align*}
g_1(\tau) &= \prod_{j=0}^1\dfrac{\eta_{40,10j+2}(\tau)\eta_{40,10j+8}(\tau)}{\eta_{40,10j+1}^2(\tau)\eta_{40,10j+4}^2(\tau)\eta_{40,10j+6}^2(\tau)\eta_{40,10j+9}^2(\tau)},\\
g_2(\tau) &= \prod_{j=0}^1\dfrac{\eta_{40,10j+4}(\tau)\eta_{40,10j+6}(\tau)}{\eta_{40,10j+2}^2(\tau)\eta_{40,10j+3}^2(\tau)\eta_{40,10j+7}^2(\tau)\eta_{40,10j+8}^2(\tau)},
\end{align*}
so by Lemma \ref{lem25}, we see that $g_1(\tau)$ and $g_2(\tau)$ are modular functions on $\Gamma_1(40)$. In view of (\ref{eq3}) and (\ref{eq4}), we deduce that $T_1(\tau)$ and $T_2(\tau)$ are modular functions on $\Gamma_1(40)$.
\end{proof}

We now show the modular properties of the functions $g(\tau)$ and $U^2(\tau)$.

\begin{theorem}\label{thm35}
The functions $g(\tau)$ and $U^2(\tau)$ both generate the field of all modular functions on $\Gamma_0(10)$ and satisfy 
\begin{equation}
g(\tau)=\dfrac{U^2(\tau)-1}{5U^2(\tau)-1}.\label{eq5}
\end{equation}
In addition, $U(2\tau)$ is a modular function on $\Gamma_0(20)$ with character $(\frac{5}{\cdot})$.
\end{theorem} 

\begin{proof}
Writing $U^2(\tau) = \eta^2(\tau)\eta^{-4}(2\tau)\eta^{-2}(5\tau)\eta^4(10\tau)$ from Lemma \ref{lem32}, we see from Lemma \ref{lem22} that $g(\tau)$ and $U^2(\tau)$ are modular functions on $\Gamma_0(10)$. By Lemma \ref{lem21}, we choose $S_{\Gamma_0(10)} = \{\infty, 0,1/2,1/5\}$ as the set of inequivalent cusps of $\Gamma_0(10)$. We compute the orders of $g(\tau)$ and $U^2(\tau)$ at each element of $S_{\Gamma_0(10)}$ using Lemma \ref{lem23}, as shown in Table \ref{tab1}.

\begin{table}[ht]
\begin{center}
	\begin{tabular}{@{}ccccc@{}}
		\hline
		cusp $r$ & $\infty$ & $0$ & $1/2$ & $1/5$\\
		\hline
		$\mbox{ord}_r g(\tau)$ & $0$ & $-1$ & $0$ & $1$ \\
		\hline
		$\mbox{ord}_r U^2(\tau)$ & $1$ & $0$ & $-1$ & $0$ \\
		\hline
	\end{tabular}
	\caption{The orders of $g(\tau)$ and $U^2(\tau)$ at the cusps of $\Gamma_0(10)$}\label{tab1}%
\end{center}
\end{table}

We see that $g(\tau)$ and $U^2(\tau)$ have simple poles at $0$ and $1/2$, respectively. Thus, we obtain $A_0(\Gamma_0(10)) = \mathbb{C}(g(\tau))=\mathbb{C}(U^2(\tau))$ and 
\begin{align*}
g(\tau) = \dfrac{aU^2(\tau)+b}{cU^2(\tau)+d}
\end{align*}
for some constants $a,b,c$ and $d$. Using the $q$-expansions
\begin{align}
g(\tau) &= 1 + 4 q + 12 q^2 + 32 q^3 + 76 q^4 + 164 q^5 + 336 q^6 + O(q^7),\label{eq6}\\
U^2(\tau) &= q-2q^2+3q^3-6q^4+11q^5-16q^6+O(q^7)\label{eq7},
\end{align}
we may take $(a,b,c,d)=(1,-1,5,-1)$, which yields (\ref{eq5}). On the other hand, by Lemma \ref{lem22}, $u(2\tau)=\eta(2\tau)\eta^{-2}(4\tau)\eta^{-1}(10\tau)\eta^2(20\tau)$ is a modular function on $\Gamma_0(20)$ with character $(\frac{5}{\cdot})$.
\end{proof}

We now derive modular functions on $\Gamma_0(10)$ involving $I(\tau), J(\tau), T_1(\tau)$ and $T_2(\tau)$. As a consequence, we provide another proof of an identity between $I(\tau)$ and $J(\tau)$, which was first established by Rajkhowa and Saikia \cite{rajsai}, and present an analogous identity between $T_1(\tau)$ and $T_2(\tau)$. We use the fact that $\Gamma_0(10)$ is generated by $\Gamma_1(10)$ and the matrix $\gamma:=[\begin{smallmatrix}
3 & -1\\10 & -3
\end{smallmatrix}]$.

\begin{theorem}\label{thm36}
The functions $I^4(\tau)+J^{-4}(\tau)$ and $J^4(\tau)+I^{-4}(\tau)$ are modular on $\Gamma_0(10)$ and satisfy 
\begin{align*}
I^4(\tau)+\dfrac{1}{J^4(\tau)} &= \dfrac{1}{U^2(\tau)}-2+3U^2(\tau),\\
J^4(\tau)+\dfrac{1}{I^4(\tau)} &= \dfrac{1}{U^6(\tau)}-\dfrac{2}{U^4(\tau)}+\dfrac{3}{U^2(\tau)}.
\end{align*}
\end{theorem}

\begin{proof}
Since $(I^4\circ \gamma)(\tau) = J^{-4}(\tau)$ by Lemma \ref{lem24}, we see that $h(\tau):=I^4(\tau)+J^{-4}(\tau)$ and $J^4(\tau)+I^{-4}(\tau)$ are both modular on $\Gamma_0(10)$. We again take $S_{\Gamma_0(10)}=\{\infty, 0,1/2,1/5\}$ as the set of inequivalent cusps of $\Gamma_0(10)$. We calculate the orders of $I^4(\tau), J^{-4}(\tau)$ and $h(\tau)$ at each element of $S_{\Gamma_0(10)}$ using Lemmata \ref{lem21} and \ref{lem26}, as shown in Table \ref{tab2}. Here, we use the fact that for each cusp $r\in S_{\Gamma_0(10)}$,
\begin{align*}
\mbox{ord}_r h(\tau) \geq \min\{\mbox{ord}_r I^4(\tau),\mbox{ord}_r J^{-4}(\tau)\}
\end{align*}
with equality holds whenever $\mbox{ord}_r I^4(\tau)\neq \mbox{ord}_r J^{-4}(\tau)$.

\begin{table}[ht]
\begin{center}
	\begin{tabular}{@{}ccccc@{}}
		\hline
		cusp $r$ & $\infty$ & $0$ & $1/2$ & $1/5$\\
		\hline
		$\mbox{ord}_r I^4(\tau)$ & $3$ & $0$ & $-1$ & $0$ \\
		\hline
		$\mbox{ord}_r J^{-4}(\tau)$ & $-1$ & $0$ & $-1$ & $0$ \\
		\hline
		$\mbox{ord}_r h(\tau)$ & $-1$ & $\geq 0$ & $\geq -1$ & $\geq 0$ \\
		\hline
	\end{tabular}
\end{center}
\caption{The orders of $I^4(\tau), J^{-4}(\tau)$ and $h(\tau)$ at the cusps of $\Gamma_0(10)$}\label{tab2}%
\end{table}

We recall from Table \ref{tab1} that $U^2(\tau)$ has a simple pole at $1/2$ and a simple zero at $\infty$. Thus, we find that $U^2(\tau)h(\tau)$, being a nonconstant function, has a pole only at $1/2$ of order at most two and holomorphic elsewhere. Thus, by \cite[Lemma 2]{yang2}, we deduce that $U^2(\tau)h(\tau)$ is a polynomial in $U^2(\tau)$ of degree at most $2$ with complex coefficients. Writing $U^2(\tau)h(\tau)=a+bU^2(\tau)+cU^4(\tau)$ for some constants $a,b$ and $c$, and using the $q$-expansions
\begin{align*}
I^4(\tau) &= q^3 - 4 q^4 + 6 q^5 - 4 q^6 + 5 q^7 - 16 q^8+O(q^9),\\
J^4(\tau) &= q - 4 q^3 + 4 q^4 + 6 q^5 - 16 q^6 + 6 q^7 + 28 q^8+O(q^9)
\end{align*}
and (\ref{eq7}), we obtain $(a,b,c)=(1,-2,3)$. Hence, we arrive at $U^2(\tau)h(\tau) = 1-2U^2(\tau)+3U^4(\tau)$, which is equivalent to the first identity. The second identity follows from dividing both sides of the first identity by $U^4(\tau)=I^4(\tau)J^{-4}(\tau)$.
\end{proof}

\begin{corollary}[\cite{rajsai}]\label{cor37} We have $I(\tau)(1+I(\tau)J(\tau))=J^2(\tau)(I^3(\tau)+J(\tau))$.
\end{corollary}

\begin{proof}
Using the definition of $U(\tau)$, we write the first identity of Theorem \ref{thm36} as
\begin{align*}
I^4(\tau)+\dfrac{1}{J^4(\tau)} = \dfrac{J^2(\tau)}{I^2(\tau)}-2+3\dfrac{I^2(\tau)}{J^2(\tau)},
\end{align*}
which, after clearing denominators, is equivalent to $F(I(\tau), J(\tau))=0$, where
\begin{align*}
F(X,Y) = (-X+X^2Y+X^3Y^2-Y^3)(-X-X^2Y+X^3Y^2+Y^3).
\end{align*}
Using the $q$-expansions of $I(\tau)$ and $J(\tau)$, we see that the first factor of $F(X,Y)$ does not vanish when $(X,Y)=(I(\tau),J(\tau))$ while the second factor does. 
We thus get $-I(\tau)-I^2(\tau)J(\tau)+I^3(\tau)J^2(\tau)+J^3(\tau)=0$, which yields the desired identity.
\end{proof}

\begin{theorem}\label{thm38}
The function $g_1^4(\tau)+g_2^4(\tau)$ is modular on $\Gamma_0(10)$ and satisfies 
\begin{equation*}
g_1^4(\tau)+g_2^4(\tau) = \dfrac{1}{g(\tau)}-2+3g(\tau).
\end{equation*}
\end{theorem}

\begin{proof}
We write the functions $g_1^4(\tau) = \eta_{10,1}^{-8}(\tau)\eta_{10,2}^4(\tau)\eta_{10,4}^{-8}(\tau)$ and $g_2^4(\tau) = \eta_{10,2}^{-8}(\tau)\eta_{10,4}^4(\tau)\eta_{10,3}^{-8}(\tau)$. Since $(g_1^4\circ \gamma)(\tau) = g_2^4(\tau)$ by Lemma \ref{lem24}, we know that $f(\tau):=g_1^4(\tau)+g_2^4(\tau)$ is a modular function on $\Gamma_0(10)$. By choosing again the set $S_{\Gamma_0(10)}=\{\infty, 0,1/2,1/5\}$ of inequivalent cusps of $\Gamma_0(10)$, we compute the orders of $g_1^4(\tau), g_2^4(\tau)$ and $f(\tau)$ at each element of $S_{\Gamma_0(10)}$ using Lemmata \ref{lem21} and \ref{lem26}, as shown in Table \ref{tab3}. We use the fact that for each cusp $r\in S_{\Gamma_0(10)}$,
\begin{align*}
\mbox{ord}_r f(\tau) \geq \min\{\mbox{ord}_r g_1^4(\tau),\mbox{ord}_r g_2^4(\tau)\}
\end{align*}
with equality holds whenever $\mbox{ord}_r g_1^4(\tau)\neq \mbox{ord}_r g_2^4(\tau)$.

\begin{table}[ht]
\begin{center}
	\begin{tabular}{@{}ccccc@{}}
		\hline
		cusp $r$ & $\infty$ & $0$ & $1/2$ & $1/5$\\
		\hline
		$\mbox{ord}_r g_1^4(\tau)$ & $0$ & $-1$ & $0$ & $-1$ \\
		\hline
		$\mbox{ord}_r g_2^4(\tau)$ & $0$ & $-1$ & $0$ & $3$ \\
		\hline
		$\mbox{ord}_r f(\tau)$ & $\geq 0$ & $\geq -1$ & $\geq 0$ & $-1$ \\
		\hline
	\end{tabular}
\end{center}
\caption{The orders of $g_1^4(\tau), g_2^4(\tau)$ and $f(\tau)$ at the cusps of $\Gamma_0(10)$}\label{tab3}%
\end{table}

We know from Table \ref{tab1} that $g(\tau)$ has a simple pole at $0$ and a simple zero at $1/5$. Thus, we see that the nonconstant function $g(\tau)f(\tau)$ has a pole only at $0$ of order at most two and holomorphic elsewhere. By \cite[Lemma 2]{yang2}, we deduce that $g(\tau)f(\tau)$ is a polynomial in $g(\tau)$ of degree at most two with complex coefficients. Setting $g(\tau)f(\tau) = a+bg(\tau)+cg^2(\tau)$ for some constants $a,b$ and $c$, and applying the $q$-expansions
\begin{align*}
g_1^4(\tau) &= 1 + 8 q + 32 q^2 + 88 q^3 + 200 q^4 + 424 q^5 + 872 q^6+O(q^7),\\
g_2^4(\tau) &= 1 + 8 q^2 + 8 q^3 + 32 q^4 + 64 q^5 + 120 q^6+O(q^7)
\end{align*}
and (\ref{eq6}), we get $(a,b,c)=(1,-2,3)$. Hence, we arrive at $g(\tau)f(\tau) = 1-2g(\tau)+3g^2(\tau)$, which is equivalent to the desired identity.
\end{proof}

\begin{corollary}\label{cor39}
We have the identity
\begin{equation*}
T_1(\tau)(T_1(\tau)+T_2(\tau))(1+T_2^2(\tau))=T_2(\tau)(1+T_1(\tau)T_2(\tau))(1+T_1^2(\tau)).
\end{equation*}
\end{corollary}

\begin{proof}
Using the equations (\ref{eq3}) and (\ref{eq4}) and the definition of $g(\tau)$, we transform the identity in Theorem \ref{thm38} into
\begin{align*}
\left(\dfrac{1+T_1(\tau)}{1-T_1(\tau)}\right)^4+\left(\dfrac{1+T_2(\tau)}{1-T_2(\tau)}\right)^4 &= 3\left(\dfrac{1+T_1(\tau)}{1-T_1(\tau)}\right)^2\left(\dfrac{1+T_2(\tau)}{1-T_2(\tau)}\right)^2 - 2\\
&+\left(\dfrac{1-T_1(\tau)}{1+T_1(\tau)}\right)^2\left(\dfrac{1-T_2(\tau)}{1+T_2(\tau)}\right)^2.
\end{align*}
Clearing denominators, the above identity is equivalent to $H(T_1(\tau),T_2(\tau))=0$, where
\begin{align*}
H(X,Y) &= (-X^2+Y-XY+X^2Y+XY^2-X^2Y^2+X^3Y^2-XY^3)\\
&\times (-X+XY-X^2Y+X^3Y+Y^2-XY^2+X^2Y^2-X^2Y^3).
\end{align*}
Applying the $q$-expansions 
\begin{align}
T_1(\tau) &= q + q^4 - q^8 - q^{11} - q^{14} + q^{15}+O(q^{16}),\label{eq8}\\
T_2(\tau) &= q^2 + q^3 - q^{11} - q^{12} - q^{13} - q^{14}+O(q^{16})\nonumber
\end{align}
we get that the first factor of $H(X,Y)$ vanishes when $(X,Y) = (T_1(\tau),T_2(\tau))$ while the second factor does not. Hence, we arrive at
\begin{align*}
-T_1^2(\tau)+T_2(\tau)&-T_1(\tau)T_2(\tau)+T_1^2(\tau)T_2(\tau)\\
&+T_1(\tau)T_2^2(\tau)-T_1^2(\tau)T_2^2(\tau)+T_1^3(\tau)T_2^2(\tau)-T_1(\tau)T_2^3(\tau)=0,
\end{align*}
which yields the desired identity.
\end{proof}

\section{Modular equations for $g(\tau)$ and $U(\tau)$}\label{sec4}

We recall in Section \ref{sec3} that the function $g(\tau)$ generates the field $A_0(\Gamma_0(10))$. The following lemma says that we can define an affine model for the modular curve $X_0(10n)$ over $\mathbb{Q}$ using $g(\tau)$. This affine model gives rise to modular equations for $g(\tau)$ and $U(\tau)$ of suitable level. We omit the proof as it is similar to that of \cite[Lemma 2.6]{leeparkd}. 

\begin{lemma}\label{lem41}
For all positive integers $n$, we have $A_0(\Gamma_0(10n))_{\mathbb{Q}}= \mathbb{Q}(g(\tau),g(n\tau))$.
\end{lemma}

We next describe the behavior of $g(\tau)$ at all cusps in $\mathbb{Q}\cup\{\infty\}$ using Table \ref{tab1}.

\begin{lemma}\label{lem42}
Let $a, c, a', c'$ and $n$ be integers with $n$ positive. Then:
\begin{enumerate}
\item[(1)] The modular function $g(\tau)$ has a pole at $a/c\in \mathbb{Q}\cup\{\infty\}$ if and only if $\gcd(a,c) = 1$ and $c\equiv \pm1, \pm 3\pmod {10}$.
\item[(2)] The modular function $g(n\tau)$ has a pole at $a'/c'\in \mathbb{Q}\cup\{\infty\}$ if and only if there exist integers $a$ and $c$ such that $a/c=na'/c', \gcd(a,c)=1$ and $c\equiv \pm1, \pm 3\pmod {10}$.
\item[(3)] The modular function $g(\tau)$ has a zero at $a/c\in \mathbb{Q}\cup\{\infty\}$ if and only if $\gcd(a,c) = 1$ and $c\equiv 5\pmod{10}$.
\item[(4)] The modular function $g(n\tau)$ has a zero at $a'/c'\in \mathbb{Q}\cup\{\infty\}$ if and only if there exist integers $a$ and $c$ such that $a/c=na'/c', \gcd(a,c)=1$ and $c\equiv 5\pmod{10}$.
\end{enumerate}
\end{lemma}

\begin{proof}
Using Table \ref{tab1}, we know that $g(\tau)$ has a simple pole (resp., zero) at $a/c\in \mathbb{Q}\cup\{\infty\}$ which is equivalent to $0$ (resp., $1/5$) on $\Gamma_0(10)$. By Lemma \ref{lem21}, $a/c$ is equivalent to $0$ on $\Gamma_0(10)$ if and only if $(a,c)\equiv (n, s)\pmod{10}$ for some integer $n$ and $s\in(\mathbb{Z}/10\mathbb{Z})^\times$. This implies that $\gcd(a,c)=1$ and $c\equiv \pm 1,\pm 3\pmod{10}$, proving (1). On the other hand, $a/c$ is equivalent to $1/5$ on $\Gamma_0(10)$ if and only if $(a,c)\equiv (s^{-1}+5n, 5s)\pmod{10}$ for some integer $n$. This implies that $\gcd(a,c)=1$ and $c\equiv 5\pmod{10}$, proving (3). Statements (2) and (4) immediately follow from (1) and (3).
\end{proof}

We now present the following result due to Ishida and Ishii \cite{ishdai} about the coefficients of the modular equation between two functions $f_1(\tau)$ and $f_2(\tau)$ that generate the field $A_0(\Gamma')$ for some congruence subgroup $\Gamma'$ of $\mbox{SL}_2(\mathbb{Z})$. This can be used to derive modular equations for $g(\tau)$ and $U(\tau)$.

\begin{proposition}[\cite{chokimkoo,ishdai}]\label{prop43}
Let $\Gamma'$ be a congruence subgroup of $\mbox{SL}_2(\mathbb{Z})$, and $f_1(\tau)$ and $f_2(\tau)$ be nonconstant functions with $A_0(\Gamma')=\mathbb{C}(f_1(\tau),f_2(\tau))$. For $k\in\{1,2\}$, let $d_k$ be the total degree of poles of $f_k(\tau)$. Let 
\begin{align*}
F(X,Y)=\sum_{\substack{0\leq i\leq d_2\\0\leq j\leq d_1}}C_{i,j}X^iY^j\in \mathbb{C}[X,Y]
\end{align*}
satisfy $F(f_1(\tau),f_2(\tau))=0$. Let $S_{\Gamma'}$ be the set of all the inequivalent cusps of $\Gamma'$, and for $k\in\{1,2\}$, define 
\begin{align*}
S_{k,0}&= \{r\in S_{\Gamma'} : f_k(\tau)\text{ has a zero at }r\},\\
S_{k,\infty}&= \{r\in S_{\Gamma'} : f_k(\tau)\text{ has a pole at }r\}.
\end{align*}
Further let
\begin{align*}
a = -\sum_{r\in S_{1,\infty}\cap S_{2,0}}\mbox{ord}_r f_1(\tau)\quad\text{ and }\quad b = \sum_{r\in S_{1,0}\cap S_{2,0}}\mbox{ord}_r f_1(\tau)
\end{align*}
with the convention that $a=0$ (resp. $b=0$) whenever $S_{1,\infty}\cap S_{2,0}$ (resp. $S_{1,0}\cap S_{2,0}$) is empty. Then
\begin{enumerate}
\item[(1)] $C_{d_2,a}\neq 0$ and if, in addition, $S_{1,\infty}\subseteq S_{2,\infty}\cup S_{2,0}$, then $C_{d_2,j}=0$ for all $j\neq a$;
\item[(2)] $C_{0,b}\neq 0$ and if, in addition, $S_{1,0}\subseteq S_{2,\infty}\cup S_{2,0}$, then $C_{0,j}=0$ for all $j\neq b$.
\end{enumerate}
By interchanging the roles of $f_1(\tau)$ and $f_2(\tau)$, one may obtain properties analogous to (1) - (2).
\end{proposition}

We now apply Lemma \ref{lem41} and Proposition \ref{prop43} to prove that modular equations for $g(\tau)$ of any level $n$ exist for positive integers $n$, which serve as affine models for $X_0(10n)$.

\begin{theorem}\label{thm44}
One can explicitly obtain modular equations for $g(\tau)$ and $U(\tau)$ of level $n$ for all positive integers $n$.
\end{theorem}

\begin{proof}
Following the proof of \cite[Theorem 1.4]{leeparkd}, we infer from Lemma \ref{lem41} that there is a polynomial $G_n(X,Y)\in\mathbb{C}[X,Y]$ such that $G_n(g(\tau),g(n\tau))=0, \deg_X G_n(X,Y)=d_n$ and $\deg_Y G_n(X,Y)=d_1$, where for $j\in\{1,n\}$, $d_j$ is the total degree of poles of $g(j\tau)$ on $\Gamma_0(10n)$ (which is also equal to the extension degree $[A_0(\Gamma_0(10n)):\mathbb{C}(g(j\tau))]$). Moreover, the coefficients of $G_n(X,Y)$ are rational numbers because $g(\tau)$ has a rational $q$-expansion (\ref{eq6}). In view of the proof of \cite[Lemma 3.1]{chokimkoo}, we see that $G_n(X,g(n\tau))$ (resp. $G_n(g(\tau),Y)$) is the minimal polynomial of $g(\tau)$ (resp. $g(n\tau)$) over $\mathbb{C}(g(n\tau))$ (resp. $\mathbb{C}(g(\tau))$) of degree $d_n$ (resp. $d_1$). Thus, we get that $G_n(g(\tau),g(n\tau))=0$ is the modular equation for $g(\tau)$ of level $n$ for all positive integers $n$, so by (\ref{eq5}), we know that 
\begin{equation*}
G_n\left(\dfrac{U^2(\tau)-1}{5U^2(\tau)-1}, \dfrac{U^2(n\tau)-1}{5U^2(n\tau)-1}\right)=0
\end{equation*}
for all positive integers $n$. Consider the polynomial 
\begin{equation*}
\mathcal{G}_n(X,Y) = (5X^2-1)^{d_n}(5Y^2-1)^{d_1}G_n\left(\dfrac{X^2-1}{5X^2-1},\dfrac{Y^2-1}{5Y^2-1}\right).
\end{equation*}
Since $\mathcal{G}_n(U(\tau),U(n\tau))=0$, we choose an irreducible factor $U_n(X,Y)\in\mathbb{Q}[X,Y]$ of $\mathcal{G}_n(X,Y)$ such that $U_n(U(\tau),U(n\tau))=0$ as follows: after substituting the $q$-expansion of $U(\tau)$ into the irreducible factors of $\mathcal{G}_n(X,Y)$, we find which among these becomes $O(q^N)$ for some positive integer $N$. We thus conclude that $U_n(U(\tau),U(n\tau))=0$ is the modular equation for $U(\tau)$ of level $n$ for all positive integers $n$.
\end{proof}

We now use Theorem \ref{thm44} to derive explicitly the modular equations for $g(\tau)$ and $U(\tau)$ of levels two and five, and provide information about the coefficients of the modular equations for $g(\tau)$ of odd prime levels $p\neq 5$.

\begin{theorem}[Modular equations of level two]\label{thm45} We have
\begin{enumerate}
\item[(1)] $g(\tau)-g(\tau)^2-4g(\tau)g(2\tau)-g(2\tau)^2+5g(\tau)g(2\tau)^2=0$
\item[(2)] $U(\tau)^4+4U(\tau)^2U(2\tau)^2+U(2\tau)^4=U(2\tau)^2(1+5U(\tau)^4)$
\end{enumerate}
\end{theorem}

\begin{proof}
In view of Lemmata \ref{lem41} and \ref{lem42}, we work on the congruence subgroup $\Gamma_0(20)$ and the following cusps: $0, 1/2, 1/5$ and $1/10$. Using Lemma \ref{lem23}, we know that $g(\tau)$ has a double pole at $0$ and a double zero at $1/5$. We also find that $g(2\tau)$ has two simple poles at $0$ and $1/2$, and two simple zeros at $1/5$ and $1/10$. Thus, we see that the total degree of poles of both $g(\tau)$ and $g(2\tau)$ are $2$, so by Proposition \ref{prop43}, there is a polynomial 
\begin{align*}
G_2(X,Y) = \sum_{0\leq i, j\leq 2} C_{i,j}X^iY^j\in\mathbb{C}[X,Y]
\end{align*}
such that $G_2(g(\tau),g(2\tau))=0$ with $C_{2,0}\neq 0$ and $C_{0,2}\neq 0$. We also get $C_{2,1}=C_{2,2}=C_{0,0}=C_{0,1}=0$. Using the $q$-expansion (\ref{eq6}) of $g(\tau)$ and taking $C_{2,0}=-1$, we arrive at $G_2(X,Y)=X-X^2-4XY-Y^2+5XY^2$ and the first identity follows from $G_2(g(\tau),g(2\tau))=0$. Using (\ref{eq5}), we see that
\begin{align*}
G_2\left(\dfrac{X^2-1}{5X^2-1},\dfrac{Y^2-1}{5Y^2-1}\right)=\dfrac{-16U_2(X,Y)}{(5X^2-1)^2(5Y^2-1)^2},
\end{align*}
where $U_2(X,Y) = X^4 - Y^2 + 4 X^2 Y^2 - 5 X^4 Y^2 + Y^4$ is irreducible. Hence, we find that $U_2(U(\tau),U(2\tau))=0$ is the modular equation for $U(\tau)$ of level two, and the second identity follows.
\end{proof}

\begin{theorem}[Modular equations of level five]\label{thm46} We have
\begin{enumerate}
\item[(1)] $g(\tau)^5-g(5\tau)-5g(\tau)g(5\tau)P(g(\tau),g(5\tau))=25g(\tau)^2g(5\tau)^2Q(g(\tau),g(5\tau))$, where
{\small \begin{align*}
		P(X,Y) &= Y+11X^3-12X^2+7X-2,\\
		Q(X,Y) &= 25X^2Y^3-5X(10X-1)Y^2+(35X^2-10X+1)Y-12X^2+7X-2.
\end{align*}}%
\item[(2)] $U(\tau)^5-U(5\tau)=5U(\tau)U(5\tau)R(U(\tau),U(5\tau))$, where
\begin{align*}
	R(X,Y) = (5Y^4-5Y^2+1)X^3-Y(5Y^2-3)X^2+(3Y^2-1)X-Y.
\end{align*}
\end{enumerate}
\end{theorem}

\begin{proof}
Invoking Lemmata \ref{lem41} and \ref{lem42}, we focus on the congruence subgroup $\Gamma_0(50)$ and the following cusps: $0, 1/5,2/5,3/5,4/5$ and $1/25$. Using Lemma \ref{lem23}, we know that $g(\tau)$ has a pole at $0$ of order five and five simple zeros at $1/5, 2/5, 3/5, 4/5$ and $1/25$. We also know that $g(5\tau)$ has a zero at $1/25$ of order five and five simple poles at $0, 1/5, 2/5, 3/5$ and $4/5$. We see that the total degree of poles of both $g(\tau)$ and $g(5\tau)$ are $5$. Thus, Proposition \ref{prop43} tells us there is a polynomial 
\begin{align*}
G_5(X,Y) = \sum_{0\leq i, j\leq 5} C_{i,j}X^iY^j\in\mathbb{C}[X,Y]
\end{align*}
such that $G_5(g(\tau),g(5\tau))=0$. Moreover, we get $C_{5,0}\neq 0, C_{0,1}\neq 0$ and $C_{5,j}=C_{0,k}=0$ for $j\in\{1,2,3,4,5\}$ and $k\in\{0,2,3,4,5\}$. Switching the roles of $g(\tau)$ and $g(5\tau)$ yields $C_{4,5}\neq 0$ and $C_{1,5}=C_{2,5}=C_{3,5}=C_{j,0}=0$ for $j\in\{1,2,3,4\}$. Using the $q$-expansion (\ref{eq6}) of $g(\tau)$ and taking $C_{5,0}=1$, we arrive at
\begin{align*}
G_5(X,Y)=X^5-Y-5XYP(X,Y)-25X^2Y^2Q(X,Y),
\end{align*}
and the first identity follows from $G_5(g(\tau),g(5\tau))=0$. Using (\ref{eq5}), we get
\begin{align*}
G_5\left(\dfrac{X^2-1}{5X^2-1},\dfrac{Y^2-1}{5Y^2-1}\right)=\dfrac{1024\mathcal{G}_5(X,Y)}{(5X^2-1)^5(5Y^2-1)^5},
\end{align*}
where $\mathcal{G}_5(X,Y)$ is the product of two irreducible polynomials
{\small \begin{align*}
	&X^5 - Y + 5 X^2 Y - 5 X^4 Y + 5 X Y^2 - 15 X^3 Y^2 - 15 X^2 Y^3 + 
	25 X^4 Y^3 + 25 X^3 Y^4 - 25 X^4 Y^5,\\
	&X^5 + Y - 5 X^2 Y + 5 X^4 Y + 5 X Y^2 - 15 X^3 Y^2 + 15 X^2 Y^3 - 
	25 X^4 Y^3 + 25 X^3 Y^4 + 25 X^4 Y^5.
\end{align*}}%
Using the $q$-expansion of $U(\tau)$, we see that the first factor of $\mathcal{G}_5(X,Y)$ vanishes when $(X,Y)=(U(\tau),U(5\tau))$, while the second factor does not. Hence, we take $U_5(X,Y)$ as the first factor of $\mathcal{G}_5(X,Y)$, so that $U_5(U(\tau),U(5\tau))=0$ is the modular equation for $U(\tau)$ of level five. The second identity now follows from this equation.
\end{proof}

\begin{theorem}\label{thm47}
Let $G_p(X,Y) = \sum_{0\leq i, j\leq p+1} C_{i,j} X^iY^j\in \mathbb{Q}[X,Y]$ be the modular equation for $g(\tau)$ of odd prime level $p\neq 5$. Then
\begin{enumerate}
\item[(1)] $C_{p+1,0}\neq 0$ and $C_{0,p+1}\neq 0$,
\item[(2)] $C_{p+1,j}  = C_{j,p+1} = 0$ for $j\in \{1,2,\ldots, p+1\}$, and
\item[(3)] $C_{0,j}  = C_{j,0} = 0$ for $j\in \{0,1,\ldots, p\}$.
\end{enumerate}
\end{theorem}

\begin{proof}
In view of Lemmata \ref{lem41} and \ref{lem42}, we work on the congruence subgroup $\Gamma_0(10p)$ and the following cusps: $0, 1/5, 1/p$ and $1/5p$. We infer from Lemma \ref{lem42} that both $g(\tau)$ and $g(p\tau)$ have poles at $0$ and $1/p$ and zeros at $1/5$ and $1/5p$. Using Lemma \ref{lem23}, we compute
\begin{align*}
\mbox{ord}_0g(\tau)=\mbox{ord}_{1/p}g(p\tau)=-p\quad\text{ and }\quad\mbox{ord}_{1/p}g(\tau)=\mbox{ord}_0g(p\tau)=-1,	
\end{align*}
so the total degrees of poles of both $g(\tau)$ and $g(p\tau)$ are $p+1$. We see from Proposition \ref{prop43} that there is a polynomial
\begin{align*}
G_p(X,Y) = \sum_{0\leq i, j\leq p+1}C_{i,j}X^iY^j\in \mathbb{Q}[X,Y]	
\end{align*}
such that $G_p(g(\tau),g(p\tau))=0$. We also have $C_{p+1,0}\neq 0$ and $C_{p+1,j}=0$ for $j\in \{1,2,\ldots, p+1\}$. Since $\mbox{ord}_{1/5}g(\tau)=p$ and $\mbox{ord}_{1/5p}g(\tau)=1$ from Lemma \ref{lem23}, we get $C_{0,p+1}\neq 0$ and $C_{0,j}=0$ for $j\in \{0,1,\ldots, p\}$. Interchanging the roles of $g(\tau)$ and $g(p\tau)$, we obtain $C_{j,p+1}=0$ for $j\in \{1,2,\ldots, p+1\}$ and $C_{j,0}=0$ for $j\in \{0,1,\ldots, p\}$.
\end{proof}

We provide the modular equations $G_p(X,Y)$ and $U_p(X,Y)$ for $g(\tau)$ and $U(\tau)$, respectively, of prime level $p\leq 11$, as shown in Tables \ref{tab4} and \ref{tab5}. We observe in Table \ref{tab4} that for such odd primes $p\neq 5$, $G_p(X,Y)$ satisfies the Kronecker congruence $G_p(X,Y)\equiv (X^p-Y)(X-Y^p)\pmod{p}$.

\begin{table}
\begin{center}
{\scriptsize \begin{tabular}{@{}ll@{}}
		\hline
		$p$ & $G_p(X,Y)$ \\
		\hline
		$2$ & $X-X^2-4XY-Y^2+5XY^2$ \\
		\hline
		$3$ & $(X^3-Y)(X-Y^3)-3XY(-X - Y + 2 X Y) (2 - 3 X - 3 Y + 4 X Y)$\\
		\hline
		$5$ & $\begin{aligned}
			&X^5-Y-5XY(-2 + 7 X - 12 X^2 + 11 X^3 + Y - 10 X Y + 35 X^2 Y - 60 X^3 Y + 5 X Y^2\\ 
			&- 50 X^2 Y^2 + 175 X^3 Y^2 + 25 X^2 Y^3 - 250 X^3 Y^3 + 125 X^3 Y^4)
		\end{aligned}$\\
		\hline
		$7$ & $\begin{aligned}
			&(X^7-Y)(X-Y^7)-7 X Y (-2 X + 11 X^2 - 32 X^3 + 56 X^4 - 56 X^5 + 24 X^6 - 2 Y\\
			&+ 32 X Y - 194 X^2 Y + 580 X^3 Y - 952 X^4 Y + 812 X^5 Y -280 X^6 Y + 11 Y^2 - 194 X Y^2\\
			&+ 1273 X^2 Y^2 - 3956 X^3 Y^2 + 6234 X^4 Y^2 - 4760 X^5 Y^2 + 1400 X^6 Y^2 - 32 Y^3\\ 
			&+ 580 X Y^3 -3956 X^2 Y^3 + 12678 X^3 Y^3 - 19780 X^4 Y^3 + 14500 X^5 Y^3 - 4000 X^6 Y^3\\
			&+ 56 Y^4 - 952 X Y^4 + 6234 X^2 Y^4 - 19780 X^3 Y^4 + 31825 X^4 Y^4 - 24250 X^5 Y^4\\
			&+ 6875 X^6 Y^4 - 56 Y^5 + 812 X Y^5 - 4760 X^2 Y^5 + 14500 X^3 Y^5 - 24250 X^4 Y^5\\
			&+ 20000 X^5 Y^5 - 6250 X^6 Y^5 + 24 Y^6 - 280 X Y^6 + 1400 X^2 Y^6 -4000 X^3 Y^6\\
			&+ 6875 X^4 Y^6 - 6250 X^5 Y^6 + 2232 X^6 Y^6)
		\end{aligned}$\\
		\hline
		$11$ & $\begin{aligned}
			&(X^{11}-Y)(X-Y^{11})-11 X Y (-2 X + 19 X^2 - 104 X^3 + 370 X^4 - 908 X^5 + 1562 X^6\\
			&-1848 X^7 + 1425 X^8 - 630 X^9 + 117 X^{10} - 2 Y + 44 X Y - 442 X^2 Y + 2696 X^3 Y\\
			&- 10860 X^4 Y + 29229 X^5 Y - 52008 X^6 Y + 59678 X^7 Y - 42690 X^8 Y + 17499 X^9 Y\\
			&- 3150 X^{10} Y + 19 Y^2 - 442 X Y^2 + 4795 X^2 Y^2 - 31432 X^3 Y^2 + 133446 X^4 Y^2\\
			&- 373080 X^5 Y^2 + 683634 X^6 Y^2 - 791680 X^7 Y^2 + 552585 X^8 Y^2 - 213450 X^9 Y^2\\
			&+ 35625 X^{10} Y^2 - 104 Y^3 + 2696 X Y^3 - 31432 X^2 Y^3 + 211646 X^3 Y^3\\
			&- 898904 X^4 Y^3 + 2521738 X^5 Y^3 - 4723040 X^6 Y^3 + 5614805 X^7 Y^3 -3958400 X^8 Y^3\\
			&+ 1491950 X^9 Y^3 - 231000 X^{10} Y^3 + 370 Y^4 - 10860 X Y^4 + 133446 X^2 Y^4\\
			&- 898904 X^3 Y^4 + 3725216 X^4 Y^4 -10256220 X^5 Y^4 + 19356124 X^6 Y^4\\
			&- 23615200 X^7 Y^4 + 17090850 X^8 Y^4 - 6501000 X^9 Y^4 + 976250 X^{10} Y^4\\
			&- 908 Y^5 + 29229 X Y^5 - 373080 X^2 Y^5 + 2521738 X^3 Y^5 - 10256220 X^4 Y^5\\
			&+ 27521182 X^5 Y^5 - 51281100 X^6 Y^5 + 63043450 X^7 Y^5 - 46635000 X^8 Y^5\\
			&+ 18268125 X^9 Y^5 - 2837500 X^{10} Y^5 + 1562 Y^6 - 52008 X Y^6 + 683634 X^2 Y^6\\
			&- 4723040 X^3 Y^6 + 19356124 X^4 Y^6 - 51281100 X^5 Y^6 + 93130400 X^6 Y^6\\
			&- 112363000 X^7 Y^6 + 83403750 X^8 Y^6 - 33937500 X^9 Y^6 + 5781250 X^{10} Y^6 - 1848 Y^7\\
			&+ 59678 X Y^7 - 791680 X^2 Y^7 + 5614805 X^3 Y^7 - 23615200 X^4 Y^7 + 63043450 X^5 Y^7\\
			&- 112363000 X^6 Y^7 + 132278750 X^7 Y^7 - 98225000 X^8 Y^7 + 42125000 X^9 Y^7\\
			&- 8125000 X^{10} Y^7 + 1425 Y^8 - 42690 X Y^8 + 552585 X^2 Y^8 - 3958400 X^3 Y^8\\
			&+ 17090850 X^4 Y^8 - 46635000 X^5 Y^8 + 83403750 X^6 Y^8 - 98225000 X^7 Y^8\\
			&+ 74921875 X^8 Y^8 - 34531250 X^9 Y^8 + 7421875 X^{10} Y^8 - 630 Y^9 + 17499 X Y^9\\
			&- 213450 X^2 Y^9 + 1491950 X^3 Y^9 - 6501000 X^4 Y^9 + 18268125 X^5 Y^9 - 33937500 X^6 Y^9\\
			&+ 42125000 X^7 Y^9 - 34531250 X^8 Y^9 + 17187500 X^9 Y^9 - 3906250 X^{10} Y^9 + 117 Y^{10}\\
			&- 3150 X Y^{10} + 35625 X^2 Y^{10} - 231000 X^3 Y^{10} + 976250 X^4 Y^{10} - 2837500 X^5 Y^{10}\\
			&+ 5781250 X^6 Y^{10} - 8125000 X^7 Y^{10} + 7421875 X^8 Y^{10} - 3906250 X^9 Y^{10} + 887784 X^{10} Y^{10})
		\end{aligned}$\\
		\hline
\end{tabular}}
\end{center}
\caption{Modular equations $G_p(X,Y)=0$ for $g(\tau)$ of prime level $p\leq 11$}\label{tab4}%
\end{table}

\begin{table}
\begin{center}
{\footnotesize \begin{tabular}{@{}ll@{}}
		\hline
		$p$ & $U_p(X,Y)$ \\
		\hline
		$2$ & $X^4-Y^2+4X^2Y^2-5X^4Y^2+Y^4$ \\
		\hline
		$3$ & $(X^3-Y)(X+Y^3)-3 X Y (-X^2 - Y^2 + 2 X^2 Y^2)$ \\
		\hline
		$5$ & $\begin{aligned}
			&X^5-Y-5XY(-X + X^3 - Y + 3 X^2 Y + 3 X Y^2 - 5 X^3 Y^2 - 5 X^2 Y^3 + 5 X^3 Y^4)
		\end{aligned}$\\
		\hline	
		$7$ & $\begin{aligned}
			&(X^7-Y)(X+Y^7)-7 X Y (-X^2 + 2 X^4 - 2 X^6 - 2 X^5 Y - Y^2 + 9 X^2 Y^2 - 20 X^4 Y^2\\
			&+ 10 X^6 Y^2 + 2 Y^4 - 20 X^2 Y^4 + 45 X^4 Y^4 - 25 X^6 Y^4 + 2 X Y^5 - 2 Y^6 + 10 X^2 Y^6\\
			&- 25 X^4 Y^6 + 18 X^6 Y^6)
		\end{aligned}$\\
		\hline
		$11$ & $\begin{aligned}
			&(X^{11}-Y)(X-Y^{11})-11 X Y (-X^2 + 4 X^4 - 8 X^6 + 9 X^8 - 3 X^{10} + 3 X^5 Y - 18 X^7 Y\\
			&+ 9 X^9 Y - Y^2 + 11 X^2 Y^2 - 56 X^4 Y^2 + 160 X^6 Y^2 - 139 X^8 Y^2 + 45 X^{10} Y^2 + 18 X^3 Y^3\\
			&- 126 X^5 Y^3 + 153 X^7 Y^3 - 90 X^9 Y^3 + 4 Y^4 - 56 X^2 Y^4 + 324 X^4 Y^4 - 800 X^6 Y^4\\
			&+ 800 X^8 Y^4 - 200 X^{10} Y^4 + 3 X Y^5 - 126 X^3 Y^5 + 594 X^5 Y^5 - 630 X^7 Y^5 + 75 X^9 Y^5\\
			&- 8 Y^6 + 160 X^2 Y^6 - 800 X^4 Y^6 + 1620 X^6 Y^6 - 1400 X^8 Y^6 + 500 X^{10} Y^6 - 18 X Y^7\\
			&+ 153 X^3 Y^7 - 630 X^5 Y^7 + 450 X^7 Y^7 + 9 Y^8 - 139 X^2 Y^8 + 800 X^4 Y^8 - 1400 X^6 Y^8\\
			&+ 1375 X^8 Y^8 - 625 X^{10} Y^8 + 9 X Y^9 - 90 X^3 Y^9 + 75 X^5 Y^9 - 3 Y^{10} + 45 X^2 Y^{10}\\
			&- 200 X^4 Y^{10} + 500 X^6 Y^{10} - 625 X^8 Y^{10} + 284 X^{10} Y^{10})
		\end{aligned}$\\
		\hline
\end{tabular}}
\end{center}
\caption{Modular equations $U_p(X,Y)=0$ for $U(\tau)$ of prime level $p\leq 11$}\label{tab5}%
\end{table}

We now study the properties of the modular equation for $g(\tau)$ of levels coprime to $10$. For an integer $a$ coprime to $10$, we define matrices $\sigma_a\in \mbox{SL}_2(\mathbb{Z})$ such that $\sigma_a\equiv [\begin{smallmatrix}
a^{-1} & 0\\ 0 & a
\end{smallmatrix}]\pmod{10}$. Then $\sigma_a\in \Gamma_0(10)$, and we may take
\begin{equation*}
\sigma_{\pm 1} =\pm\begin{bmatrix}
1 & 0\\ 0 & 1
\end{bmatrix}\quad\text{ and }\quad\sigma_{\pm 3}=\pm\begin{bmatrix}
-3 & -10\\10 & 33
\end{bmatrix}.
\end{equation*}

We note that we have the disjoint union 
\begin{align*}
\Gamma_0(10)\begin{bmatrix}
1 & 0\\0 & n
\end{bmatrix}\Gamma_0(10) = \bigsqcup_{0<a\mid n}\bigsqcup_{\substack{0\leq b < n/a\\\gcd(a,b,n/a)=1}}\Gamma_0(10)\sigma_a\begin{bmatrix}
a & b\\0 & \frac{n}{a}
\end{bmatrix}
\end{align*}
when $n$ is coprime to $10$, and $\Gamma_0(10)\backslash \Gamma_0(10)[\begin{smallmatrix}
1 & 0\\0 & n
\end{smallmatrix}]\Gamma_0(10)$ has $\psi(n):=n\prod_{p\mid n}(1+1/p)$ elements by \cite[Proposition 3.36]{shimura}. We now consider the polynomial
\begin{align*}
\Phi_n(X,\tau):= \prod_{0<a\mid n}\prod_{\substack{0\leq b < n/a\\\gcd(a,b,n/a)=1}}(X-(g\circ\alpha_{a,b})(\tau))
\end{align*}
where $\alpha_{a,b}:=\sigma_a[\begin{smallmatrix}
a & b\\0 & n/a
\end{smallmatrix}]$. Then the coefficients of $\Phi_n(X,\tau)$ are elementary symmetric functions of $g\circ\alpha_{a,b}$, so these are invariant under the action of $\Gamma_0(10)$ and thus are in $A_0(\Gamma_0(10))=\mathbb{C}(g(\tau))$. Thus, we may write $\Phi_n(X,\tau)$ as $\Phi_n(X,g(\tau))\in\mathbb{C}(g(\tau))[X]$. The following result gives the properties of the polynomial $\Phi_n(X,g(\tau))$ when $\gcd(n,10)=1$. We omit the proof as it is analogous to that of \cite[Theorem 2.9]{leeparkb}. We note that $\Phi_n(g(\frac{\tau}{n}),g(\tau))=0$, which can be seen as the modular equation for $g(\tau)$ of level $n$.

\begin{theorem}\label{thm48}
Suppose $n$ is a positive integer coprime to $10$ and $\Phi_n(X,Y)$ be the polynomial as above. Then
\begin{enumerate}
\item[(1)] $\Phi_n(X,Y)\in \mathbb{Z}[X,Y]$ and $\deg_X\Phi_n(X,Y)=\deg_Y\Phi_n(X,Y)=\psi(n)$.
\item[(2)] $\Phi_n(X,Y)$ is irreducible both as a polynomial in $X$ over $\mathbb{C}(Y)$ and as a polynomial in $Y$ over $\mathbb{C}(X)$.
\item[(3)] $\Phi_n(X,Y)=\Phi_n(Y,X)$.
\item[(4)] (Kronecker congruence) If $p\neq 5$ is an odd prime, then 
\begin{align*}
	\Phi_p(X,Y)\equiv (X^p-Y)(X-Y^p)\pmod{p\mathbb{Z}[X,Y]}.
\end{align*}
\end{enumerate}
\end{theorem}

We recall from Theorem \ref{thm35} that $U(2\tau)$ is a modular function on $\Gamma_0(20)$ with character $(\frac{5}{\cdot})$, whose corresponding modular curve has genus one, so we cannot apply the methods presented in \cite{leeparkb,leeparkd} to derive modular equations for $U(\tau)$. Nonetheless, by looking at the modular equations for $U(\tau)$ as shown in Table \ref{tab5}, we propose the following conjecture. 

\begin{conjecture}\label{conj49}
For an odd prime $p\neq 5$, the modular equation $U_p(X,Y)=0$ for $U(\tau)$ of level $p$ satisfies
\begin{equation*}
U_p(X,Y)\equiv\begin{cases}
	(X^p-Y)(X-Y^p)\pmod{p\mathbb{Z}[X,Y]}, &\text{ if }p\equiv \pm 1\pmod{10},\\
	(X^p-Y)(X+Y^p)\pmod{p\mathbb{Z}[X,Y]},&\text{ if }p\equiv \pm 3\pmod{10}.
\end{cases}
\end{equation*}
\end{conjecture}

\section{Evaluation of continued fractions of order $10$}\label{sec5}

We investigate in this section certain arithmetic properties of $g(\tau)$ and $U(\tau)$ at imaginary quadratic points $\tau\in\mathbb{H}$, which will be useful in evaluating continued fractions of order $10$. We start with the following result of Cho and Koo \cite{chokoo}, which asserts that a singular value of a modular function on $\Gamma^0(M)\cap \Gamma_0(N)$ with rational $q$-expansion, where 
\begin{equation*}
\Gamma^0(M) := \left\lbrace\begin{bmatrix}
a & b\\ c& d
\end{bmatrix} \in \mbox{SL}_2(\mathbb{Z}) : b\equiv 0\pmod M\right\rbrace,
\end{equation*}
generates the ring class field of some imaginary quadratic order. 

\begin{lemma}[\cite{chokoo}]\label{lem51}
Let $K$ be an imaginary quadratic field with discriminant $d_K$ and $\tau\in K\cap\mathbb{H}$ be a root of a primitive equation $aX^2+bX+c=0$ such that $b^2-4ac=f^2d_K$ and $a,b,c,f\in\mathbb{Z}$ with $f > 0$. Let $\mathfrak{F}$ be the field of all automorphic functions for $\Gamma^0(M)\cap \Gamma_0(N)$ with rational $q$-expansion. Then $K\cdot\mathfrak{F}(\tau)$ is the ring class field of the order in $K$ of conductor 
\[\dfrac{M/\gcd(M,c)\cdot N/\gcd(N,a)}{\gcd(M/\gcd(M,c), N/\gcd(N,a))}\cdot f.\]
\end{lemma}

We apply Lemma \ref{lem51} to prove the following result about generating the Hilbert class field over an imaginary quadratic field using $g(\tau)$.

\begin{theorem}\label{thm52}
Let $K$ be an imaginary quadratic field with discriminant $d_K$ and let $\tau\in K \cap\mathbb{H}$ be a root of the primitive polynomial $aX^2+bX+c\in\mathbb{Z}[X]$ such that $10\mid a$ and $b^2-4ac=d_K$. Then $K(g(\tau))=K(U^2(\tau))$ is the Hilbert class field of $K$.
\end{theorem}

\begin{proof}
We know from Theorem \ref{thm35} that $g(\tau)$ generates $A_0(\Gamma_0(10))_{\mathbb{Q}}$ since it has a rational $q$-expansion (\ref{eq6}). Thus, setting $M=f=1$ and $N=10$ in Lemma \ref{lem51} and using the divisibility condition on $a$, we see that $K(g(\tau))$ is the ring class field of the order in $K$ of conductor $Nf/\gcd(N,a)=1$, whence the maximal order in $K$. Hence, we conclude that $g(\tau)$ generates the Hilbert class field of $K$, and the field equality follows from (\ref{eq5}).
\end{proof}

We remark that the conditions of Theorem \ref{thm52} imply that
\begin{align}
d_K\equiv 0,1,4,9,16,20,24,25,36\pmod{40}. \label{eq9}
\end{align}
Conversely, if $d_K$ satisfies (\ref{eq9}), then we can find integers $a, b$ and $c$ with $10\mid a$ and $\gcd(a,b,c)=1$ such that the discriminant of $aX^2+bX+c$ is $d_K$. Indeed, we simply take $k$ to be the smallest integer such that $40k+d_K$ is a positive perfect square, and set $(a,b,c)=(10k,\sqrt{40k+d_K},1)$. We infer from Theorem \ref{thm52} that for such $k$ and a root $\tau\in K\cap \mathbb{H}$ of $10kX^2+\sqrt{40k+d_K}X+1=0$, $g(\tau)$ generates the Hilbert class field of $K$.

We next establish the following result about the singular values of $1/g(\tau)$ and $1/U(\tau)$ at imaginary quadratic points $\tau\in \mathbb{H}$. 

\begin{theorem}\label{thm53}
Let $K$ be an imaginary quadratic field and $\tau\in K\cap\mathbb{H}$. Then $1/g(\tau)$ and $1/U(\tau)$ are algebraic integers.
\end{theorem}

\begin{proof}
We first note that the function $f(\tau) := \eta^6(\tau)\eta^{-6}(5\tau)$ generates the field of all modular functions on $\Gamma_0(5)$ by Lemmata \ref{lem22} and \ref{lem23}, and satisfies the identity \cite[Theorem 5.26]{cooper}
\begin{align}
j(\tau) = \dfrac{(f^2(\tau)+250f(\tau)+3125)^3}{f^5(\tau)},\label{eq10}
\end{align}
where $j(\tau)$ is the modular $j$-invariant. We consider $f(2\tau)$ as a modular function on $\Gamma_0(10)$. By Lemma \ref{lem23}, we find that $f(2\tau)$ has a simple pole at $1/5$ and a double pole at $\infty$. Since $A_0(\Gamma_0(10))=\mathbb{C}(f(2\tau),g(\tau))$ by Theorem \ref{thm35}, we infer from Proposition \ref{prop43} that there is a polynomial 
\begin{align*}
F(X,Y) = \sum_{\substack{0\leq i\leq 1\\0\leq j\leq 3}}C_{i,j}X^iY^j\in\mathbb{C}[X,Y]
\end{align*}
such that $F(f(2\tau),g(\tau))=0$. Plugging in the $q$-expansions of $f(2\tau)$ and $g(\tau)$, we get $F(X,Y)=-1 + 10 Y + X Y - 25 Y^2 - 2 X Y^2 + X Y^3$ and 
\begin{align}
f(2\tau) = \dfrac{g_0(g_0-5)^2}{(g_0-1)^2},\label{eq11}
\end{align}
where $g_0 := g_0(\tau) := 1/g(\tau)$. Combining (\ref{eq10}) and (\ref{eq11}) leads to 
\begin{align}
j(2\tau)= \dfrac{(g_0^6+230g_0^5+275g_0^4-1500g_0^3+4375g_0^2-6250g_0+3125)^3}{g_0^5(g_0-1)^2(g_0-5)^{10}}.\label{eq12}
\end{align}
On the other hand, we deduce from (\ref{eq5}) that
\begin{align}
g_0 = \dfrac{5-U_0^2}{1-U_0^2},\label{eq13}
\end{align}
where $U_0 := U_0(\tau) := 1/U(\tau)$. Substituting (\ref{eq13}) into (\ref{eq12}) yields
\begin{align}
j(2\tau)= \dfrac{(U_0^{12}-10U_0^{10}+275U_0^8-1500U_0^6+4375U_0^4-6250U_0^2+3125)^3}{U_0^{20}(U_0^2-1)(U_0^2-5)^5}.\label{eq14}
\end{align}
Because the $j$-invariant $j(2\tau)$ is always an algebraic integer whenever $\tau\in K\cap\mathbb{H}$, we infer from (\ref{eq12}) and (\ref{eq14}) that $g_0(\tau)=1/g(\tau)$ and $U_0(\tau)=1/U(\tau)$ are algebraic integers.
\end{proof}

The next result is concerned about the singular values of $g(r\tau)$ and $U(r\tau)$ for positive rational numbers $r$, which can be expressed as radicals using the values of $g(\tau)$ and $U(\tau)$. We omit the proof as it is analogous to that of \cite[Theorem 1.5]{leeparkd} and \cite[Theorem 4.6]{guad}.

\begin{theorem}\label{thm54}
Given a positive rational number $r$ and a value of $g(\tau)$ expressed in terms of radicals, there exists an algorithm that expresses $g(r\tau)$ in terms of radicals, and the result also holds for $U(\tau)$. 
\end{theorem}

We describe the singular values of $J(\tau)$ and $1/I(\tau)$ at imaginary quadratic points $\tau\in\mathbb{H}$.

\begin{theorem}\label{thm55}
Let $K$ be an imaginary quadratic field and $\tau\in K\cap\mathbb{H}$. Then $J(\tau)$ and $1/I(\tau)$ are algebraic integers.
\end{theorem}

\begin{proof}
We infer from Theorem \ref{thm36} that $J^4(\tau)$ and $1/I^4(\tau)$ are the roots of the quadratic equation
\begin{align*}
T^2 - \left(\dfrac{1}{U^6(\tau)}-\dfrac{2}{U^4(\tau)}+\dfrac{3}{U^2(\tau)}\right)T + \dfrac{1}{U^4(\tau)} = 0.
\end{align*}
Owing to Theorem \ref{thm53}, we see that $J^4(\tau)$ and $1/I^4(\tau)$ are algebraic integers, which follows the desired conclusion.
\end{proof}

We now examine the singular values of $T_1(\tau)$ and $T_2(\tau)$ at imaginary quadratic points $\tau\in\mathbb{H}$. We first prove the following lemma.

\begin{lemma}\label{lem56}
Let $K$ be an imaginary quadratic field and $\tau\in K\cap\mathbb{H}$, and let $g_0(\tau) := 1/g(\tau)$. Then $4/(g_0(\tau)-1)$ and $2(g_0(\tau)+1)/(g_0(\tau)-1)$ are algebraic integers.
\end{lemma}

\begin{proof}
Let $v_0:=v_0(\tau)=4/(g_0(\tau)-1)$ and $w_0:=w_0(\tau)=2(g_0(\tau)+1)/(g_0(\tau)-1)$, so that
\begin{align*}
g_0 := g_0(\tau)= \dfrac{v_0+4}{v_0} = \dfrac{w_0+2}{w_0-2}.
\end{align*}
Substituting these values in (\ref{eq12}) yields
\begin{align*}
j(2\tau) &= \dfrac{v_0^6+4v_0^5+240v_0^4+480v_0^3+1440v_0^2+944v_0+16}{v_0(v_0+4)^5(v_0-1)^{10}}\\
&=\dfrac{w_0^6-8w_0^5+260w_0^4-1440w_0^3+4240w_0^2-6608w_0+3824}{(w_0-2)(w_0+2)^5(w_0-3)^{10}}.
\end{align*}
Since $j(2\tau)$ is an algebraic integer for $\tau\in K\cap\mathbb{H}$, we conclude that $v_0(\tau)$ and $w_0(\tau)$ are algebraic integers as well.
\end{proof}

\begin{theorem}\label{thm57}
Let $K$ be an imaginary quadratic field and $\tau\in K\cap\mathbb{H}$. Then $T_1(\tau)$ and $T_2(\tau)$ are algebraic units.
\end{theorem}

\begin{proof}
From Theorem \ref{thm38} and the definition of $g(\tau)$, we deduce that
\begin{align*}
\dfrac{1}{g_1^4(\tau)}+\dfrac{1}{g_2^4(\tau)} = g_0^3-2g_0^2+3g_0,
\end{align*}
where $g_0:=g_0(\tau)=1/g(\tau)$. Thus, the values $1/g_1(\tau)$ and $1/g_2(\tau)$ are the roots of the equation
\begin{align}\label{eq15}
T^4 - (g_0^3-2g_0^2+3g_0)T^2 + g_0^2 = 0, 
\end{align}
so by Theorem \ref{thm53}, $1/g_1(\tau)$ and $1/g_2(\tau)$ are algebraic integers. Using (\ref{eq3}) and the $q$-expansions (\ref{eq8}) of $T_1(\tau)$ and $g_1(\tau)$, we have that
\begin{align*}
\alpha:=T_1(\tau)+\dfrac{1}{T_1(\tau)}=\dfrac{2(g_1(\tau)+1)}{g_1(\tau)-1}.
\end{align*}
Replacing $T$ with $(T-2)/(T+2)$ in (\ref{eq15}), we see that $\alpha$ is the root of the equation
\begin{align*}
\left(\dfrac{T-2}{T+2}\right)^4 - (g_0^3-2g_0^2+3g_0)\left(\dfrac{T-2}{T+2}\right)^2 + g_0^2 = 0,
\end{align*}
which, after clearing denominators, is equivalent to 
\begin{align}\label{eq16}
(g_0-1)^3T^4-8(g_0^2-1)T^3-8(g_0+1)(g_0^2+3)T^2-32(g_0^2-1)T+16(g_0-1)^3=0.
\end{align}
We deduce from Lemma \ref{lem56} that 
\begin{align*}
\dfrac{8(g_0^2-1)}{(g_0-1)^3} &= \dfrac{8(g_0+1)}{(g_0-1)^2}=\dfrac{4}{g_0-1}\cdot\dfrac{2(g_0+1)}{g_0-1},\\
\dfrac{8(g_0+1)(g_0^2+3)}{(g_0-1)^3} &= \dfrac{2(g_0+1)}{g_0-1}\left(4+\dfrac{8(g_0+1)}{(g_0-1)^2}\right)
\end{align*}
are algebraic integers. Thus, in view of (\ref{eq16}), $\alpha$ actually satisfies the equation 
\begin{align*}
T^4- \dfrac{8(g_0+1)}{(g_0-1)^2}T^3-\dfrac{8(g_0+1)(g_0^2+3)}{(g_0-1)^3}T^2- \dfrac{32(g_0+1)}{(g_0-1)^2}T+16=0,
\end{align*}
which shows that $\alpha$ is an algebraic integer. As $T_1(\tau)$ and $1/T_1(\tau)$ are the roots of the equation $T^2-\alpha T+1=0$, we conclude that $T_1(\tau)$ is an algebraic unit. By similar arguments, we also show that $T_2(\tau)$ is an algebraic unit as well.
\end{proof}

Given an imaginary quadratic field $K$ with discriminant $d_K$ satisfying (\ref{eq9}) and ring of integers $\mathcal{O}_K=\mathbb{Z}[\theta]$ for some $\theta\in K\cap\mathbb{H}$, we now introduce a systematic way, based on the Shimura reciprocity law due to Cho and Koo \cite{chokoo} and Gee \cite{gee}, of computing the exact value of $g_0(\tau):=1/g(\tau)$ for some $\tau\in K\cap \mathbb{H}$, expressed in terms of $\theta$, that satisfies the conditions of Theorem \ref{thm52} by finding its minimal polynomial over $K$. We refer the interested reader to \cite{guad} for more details. Let $\mathcal{F}_N$ be the field of modular functions of level $N$ whose coefficients in their $q$-expansions lie in the cyclotomic field $\mathbb{Q}(\zeta_N)$. The main theorem of complex multiplication says that if $h\in\mathcal{F}_N$ with $h(\theta)$ finite, then $h(\theta)$ generates the ray class field $K_{(N)}$ modulo $N$ over $K$.
Consider the map $g_\theta : (\mathcal{O}_K/N\mathcal{O}_K)^\times\rightarrow \mbox{GL}_2(\mathbb{Z}/N\mathbb{Z})$ given by
\begin{equation*}
g_\theta(s\theta+t) = \begin{bmatrix}
	t - Bs & -Cs\\ s & t
\end{bmatrix}\in \mbox{GL}_2(\mathbb{Z}/N\mathbb{Z}),
\end{equation*}
where $s\theta+t\in (\mathcal{O}_K/N\mathcal{O}_K)^\times$ and $X^2+BX+C\in\mathbb{Z}[X]$ is the minimal polynomial of $\theta$ over $\mathbb{Q}$. Let $W_{N,\theta}$ be the image of $(\mathcal{O}_K/N\mathcal{O}_K)^\times$ under $g_\theta$. The Shimura reciprocity law states that there is a surjective homomorphism from $W_{N,\theta}$ to $\mbox{Gal}(K_{(N)}/K(j(\theta)))$ that sends $\alpha$ to the map $h(\theta)\mapsto h^{\alpha^{-1}}(\theta)$, whose kernel $T$ is the image $g_\theta(\mathcal{O}_K^\times)$, namely
\begin{equation*}
T = \begin{cases}
	\left\{\pm\begin{bmatrix}
		1 & 0\\ 0 & 1
	\end{bmatrix}, \pm\begin{bmatrix}
		0 &-1\\ 1 & 0
	\end{bmatrix}\right\} & \text{ if } K=\mathbb{Q}(\sqrt{-1}),\\
	\noalign{\vskip6pt}
	\left\{\pm\begin{bmatrix}
		1 & 0\\ 0 & 1
	\end{bmatrix}, \pm\begin{bmatrix}
		-1 &-1\\ 1 & 0
	\end{bmatrix}, \pm\begin{bmatrix}
		0 &-1\\ 1 & 1
	\end{bmatrix}\right\} & \text{ if } K=\mathbb{Q}(\sqrt{-3}),\\
	\noalign{\vskip6pt}
	\left\{\pm\begin{bmatrix}
		1 & 0\\ 0 & 1
	\end{bmatrix}\right\} & \text{ otherwise.}
\end{cases}
\end{equation*}

It may happen that, under certain conditions, $h(\theta)$ lies in a smaller subfield of $K_{(N)}$, as following result shows. 

\begin{lemma}[\cite{gee}]\label{lem58}
Using the above notations, suppose in addition that $\mathbb{Q}(j)\subset\mathbb{Q}(h)$. Then $h(\theta)$ generates the Hilbert class field of $K$ if and only if $W_{N,\theta}/T$ acts trivially on $h$.
\end{lemma}

It is clear from the definition of $\mathcal{O}_K$ that $\theta$ does not satisfy the conditions of Theorem \ref{thm52}. To find a suitable element in $K\cap\mathbb{H}$ that can be expressed in terms of $\theta$ and satifies the condition of that theorem, we observe that if $[\begin{smallmatrix}
a & b\\ c & d
\end{smallmatrix}]\in\mbox{SL}_2(\mathbb{Z})$ and $\beta\in K\cap\mathbb{H}$ is a root of a primitive polynomial $AX^2+BX+C\in\mathbb{Z}[X]$ with discriminant $d_K$, then $(a\beta+b)/(c\beta+d)$ is a root of a primitive polynomial 
\begin{align*}
(Ad^2-Bcd+Cc^2)X^2+(-2Abd+B(ad+bc)-2Cac)X+(Ab^2-Bab+Ca^2)
\end{align*}
with discriminant $d_K$. In particular, if $P(X)\in\mathbb{Z}[X]$ is the minimal polynomial of $\theta$ over $\mathbb{Q}$, where
\begin{align*}
P(X)=\begin{cases}
	X^2 - \dfrac{d_K}{4}, & d_K\equiv 0\pmod{4},\\
	X^2 + X+ \dfrac{1-d_K}{4}, & d_K\equiv 1\pmod{4},
\end{cases}
\end{align*}
then for an integer $k$, $-1/(\theta-k)\in K\cap \mathbb{H}$ is a root of a polynomial $P(k)X^2-P'(k)X+1$ with discriminant $d_K$. By choosing $k$ such that $10\mid P(k)$ (which exists due to (\ref{eq9})), we know from Theorem \ref{thm52} that $g_0(-1/(\theta-k))$ generates the Hilbert class field of $K$.

To get the conjugates of $g_0(-1/(\theta-k))\in K(j(\theta))$ over $K$, we now apply the action of the form class group $C(d_K)$ on $g_0(-1/(\theta-k))$ as follows. Let $x:=[a,b,c]\in C(d_K)$ be a reduced primitive binary quadratic form and $\tau_x := (-b+\sqrt{d_K})/(2a)\in\mathbb{H}$. We choose a matrix $u_x\in\mbox{GL}_2(\mathbb{Z}/10\mathbb{Z})$ such that if $d_K\equiv 0\pmod{4}$, then
\begin{equation*}
u_x \equiv\begin{cases}
	\begin{bmatrix}
		a & b/2\\ 0 & 1
	\end{bmatrix} \bmod{p} & \text{ if } p \nmid a,\\
	\noalign{\vskip6pt}
	\begin{bmatrix}
		-b/2 & -c\\ 1 & 0
	\end{bmatrix} \bmod{p} & \text{ if } p \mid a\text{ and } p\nmid c,\\
	\noalign{\vskip6pt}
	\begin{bmatrix}
		-b/2-a &-b/2-c\\ 1 & -1
	\end{bmatrix} \bmod{p} & \text{ if } p \mid a\text{ and } p\mid c
\end{cases}
\end{equation*}
and if $d_K\equiv 1\pmod{4}$, then
\begin{equation*}
u_x \equiv\begin{cases}
	\begin{bmatrix}
		a & (b-1)/2\\ 0 & 1
	\end{bmatrix} \bmod{p} & \text{ if } p \nmid a,\\
	\noalign{\vskip6pt}
	\begin{bmatrix}
		-(b+1)/2 & -c\\ 1 & 0
	\end{bmatrix} \bmod{p} & \text{ if } p \mid a\text{ and } p\nmid c,\\
	\noalign{\vskip6pt}
	\begin{bmatrix}
		-(b+1)/2-a &(1-b)/2-c\\ 1 & -1
	\end{bmatrix} \bmod{p} & \text{ if } p \mid a\text{ and } p\mid c
\end{cases}
\end{equation*}
for each prime $p$ dividing $10$. We decompose $u_x$ as $u_x = [\begin{smallmatrix}
1 & 0\\ 0 & \det u_x
\end{smallmatrix}]v_x$ with $v_x\in\mbox{SL}_2(\mathbb{Z}/10\mathbb{Z})$ and find a lift $V_x\in \mbox{SL}_2(\mathbb{Z})$ of $v_x$. We note that the action of the matrix $[\begin{smallmatrix}
1 & 0\\ 0 & \det u_x
\end{smallmatrix}]$ on $g(\tau)$ is trivial because $g(\tau)$ has a rational $q$-expansion (\ref{eq6}). We finally choose an integer $k_x$ so that the leading coefficient of the minimal polynomial of $-1/(\theta_x-k_x)$, where $\theta_x := V_x\tau_x$, over $\mathbb{Q}$ is divisible by $10$. This ensures that $g_0(-1/(\theta_x-k_x))\in K(j(\theta))$ by Theorem \ref{thm52}. Thus, by \cite[Theorem 20]{gee}, we have the action $g_0(-1/(\theta-k))^{x^{-1}}=g_0(-1/(\theta_x-k_x))$ and there is an isomorphism from $C(d_K)$ to $\mbox{Gal}(K(j(\theta))/K)$ sending the class $x^{-1}$ to the map $g_0(-1/(\theta-k))\mapsto g_0(-1/(\theta_x-k_x))$.  An application of \cite[Lemma 6.1]{chokoo} yields the following result.

\begin{proposition}\label{prop59}
Using the above notations, the set $\{g_0(-1/(\theta_x-k_x)): x\in C(d_K)\}$ is the set of all conjugates of $g_0(-1/(\theta-k))\in K(j(\theta))$ over $K$.
\end{proposition}

In view of Theorem \ref{thm53} and Proposition \ref{prop59}, we now define the minimal polynomial of $g_0(-1/(\theta-k))\in K(j(\theta))$ over $K$ given by 
\begin{align*}
H(X) := \prod_{x\in C(d_K)}(X-g_0(-1/(\theta_x-k_x)))\in \mathcal{O}_K[X].
\end{align*}

We then numerically obtain $H(X)$ using an accurate approximation of $g_0(-1/(\theta-k))$ via (\ref{eq6}).  Consequently, we get the exact value of $g_0(-1/(\theta-k))$ by choosing the correct root of $H(X)=0$, which can be used to find the values of the continued fractions of order $10$. We illustrate this with the following example.

\begin{example}\label{ex1}
Let $\theta=\sqrt{-10}$, and consider the field $K=\mathbb{Q}(\theta)$ with $d_K=-40$ and $\mathcal{O}_K=\mathbb{Z}[\theta]$. Then $g_0(-1/\theta)$ generates the Hilbert class field of $K$ by Theorem \ref{thm52}. There are two reduced primitive binary quadratic forms of discriminant $-40$, namely $[1,0,10]$ and $[2,0,5]$. We note that the action of $[1,0,10]$ on $g_0(-1/\theta)$ is trivial, so we focus on the action of $x:=[2,0,5]$ on $g_0(-1/\theta)$. We compute the matrix 
\begin{align*}
	u_x = \begin{bmatrix}
		2 & 5\\ 5 & 6
	\end{bmatrix} = \begin{bmatrix}
		1& 0\\0 & 7
	\end{bmatrix}\begin{bmatrix}
		2& 5\\5 & 8
	\end{bmatrix}\in \mbox{GL}_2(\mathbb{Z}/10\mathbb{Z}),
\end{align*}
so we take $V_x = [\begin{smallmatrix}
	22 & 35\\5 & 8
\end{smallmatrix}]$. Setting $\tau_x = \sqrt{-5/2}$, we find that the minimal polynomial of $\theta_x:=V_x\tau_x$ over $\mathbb{Q}$ is $253X^2-2220X+4870$. Thus, we deduce that $g_0(-1/\theta)^{x^{-1}}=g_0(-1/\theta_x)$, and appealing to Proposition \ref{prop59}, we compute
\begin{align*}
	H(X) = (X-g_0(-1/\theta))(X-g_0(-1/\theta_x)) \approx X^2 - 10X+5.
\end{align*}
Using the $q$-expansion (\ref{eq6}) of $g(\tau)$, we have $g_0(-1/\theta)\approx 0.5278640$, so that
\begin{align*}
	g_0(-1/\theta) = 5-2\sqrt{5},
\end{align*}
and from (\ref{eq13}), we have $U_0(-1/\theta)=\sqrt{5+2\sqrt{5}}$, where $U_0(\tau):=1/U(\tau)$. In view of the proof of Theorem \ref{thm55}, we see that $J(-1/\theta)$ and $1/I(-1/\theta)$ are the roots of the equation
\begin{align*}
	T^8-(350+156\sqrt{5})T^4+45+20\sqrt{5}=0.
\end{align*}
From the $q$-expansions of $I(\tau)$ and $J(\tau)$, we approximate $J(-1/\theta)\approx 0.5986205$ and $1/I(-1/\theta)\approx 5.1412960$, which gives 
\begin{align*}
	J(-1/\theta) &= \left(175+78\sqrt{5}-2\sqrt{15250+6820\sqrt{5}}\right)^{1/4},\\
	I(-1/\theta) &= \left(175+78\sqrt{5}+2\sqrt{15250+6820\sqrt{5}}\right)^{-1/4}.
\end{align*}
We see from (\ref{eq15}) that $1/g_1(-1/\theta)$ and $1/g_2(-1/\theta)$ satisfy the equation
\begin{align*}
	T^4-(350-156\sqrt{5})T^2+45-20\sqrt{5}=0.
\end{align*}
With $1/g_1(-1/\theta)\approx 0.5749991$ and $1/g_2(-1/\theta)\approx 0.9180259$ via the $q$-expansions of $g_1(\tau)$ and $g_2(\tau)$, we obtain 
\begin{align*}
	g_1(-1/\theta) &= \left(175-78\sqrt{5}-2\sqrt{15250-6820\sqrt{5}}\right)^{-1/2},\\
	g_2(-1/\theta) &= \left(175-78\sqrt{5}+2\sqrt{15250-6820\sqrt{5}}\right)^{-1/2},
\end{align*}
and immediately deduce the values of $T_1(-1/\theta)$ and $T_2(-1/\theta)$ using (\ref{eq3}) and (\ref{eq4}). 
On the other hand, using the modular equation for $g(\tau)$ of level two from Theorem \ref{thm45}, we verify that
\begin{align*}
	g_0(-1/(2\theta)) = (3-2\sqrt{2})(5-2\sqrt{5}).
\end{align*} 
Consequently, by performing similar calculations, we also verify that
{\footnotesize \begin{align*}
		U_0(-1/(2\theta)) &= \sqrt{\sqrt{5}+\sqrt{10}},\\
		J(-1/(2\theta)) &= \left(-15-10\sqrt{2}+19\sqrt{5}+14\sqrt{10}-2\sqrt{5(835+590\sqrt{2}-226\sqrt{5}-160\sqrt{10})}\right)^{1/4},\\
		I(-1/(2\theta)) &= \left(-15-10\sqrt{2}+19\sqrt{5}+14\sqrt{10}+2\sqrt{5(835+590\sqrt{2}-226\sqrt{5}-160\sqrt{10})}\right)^{-1/4},\\
		g_1(-1/(2\theta)) &= \left(20295 - 14350\sqrt{2} - 9074\sqrt{5} + 6416\sqrt{10}-2\sqrt{\lambda}\right)^{-1/2},\\
		g_2(-1/(2\theta)) &= \left(20295 - 14350\sqrt{2} - 9074\sqrt{5} + 6416\sqrt{10}+2\sqrt{\lambda}\right)^{-1/2},\\
\end{align*}}%
and the values of $T_1(-1/(2\theta))$ and $T_2(-1/(2\theta))$ can be obtained using (\ref{eq3}) and (\ref{eq4}), where 
\begin{align*}
	\lambda := 10(41176730 - 29116345\sqrt{2} - 18414793\sqrt{5} + 13021225\sqrt{10}).
\end{align*}
\end{example}

As shown in Example \ref{ex1}, we deduce from Theorem \ref{thm54} that for a given positive rational number $r$ and an imaginary quadratic point $\tau\in\mathbb{H}$, the values of $I(r\tau), J(r\tau), T_1(r\tau)$ and $T_2(r\tau)$ can be written in terms of radicals.

\section*{Declaration of interests}

The authors declare that they have no competing interests.

\section*{Acknowledgements}

The authors would like to thank the anonymous referees for their helpful comments that improved the presentation of this paper. 

\bibliographystyle{amsplain} 
\bibliography{contfracten}
\end{document}